\documentclass[10pt,a4paper]{article}
\usepackage[english]{babel}
\usepackage{amssymb}
\usepackage{amsmath}
\usepackage{hyperref}

\usepackage{amsthm}

\usepackage{color}

\allowdisplaybreaks[1] 

\newtheorem{prop}{Proposition}[section]
\newtheorem{lemma}[prop]{Lemma}
\newtheorem{mydef}[prop]{Definition}
\newtheorem{corol}[prop]{Corollary}
\newtheorem{mythm}[prop]{Theorem}
\newtheorem{remark}[prop]{Remark}
\newtheorem{hyp}[prop]{Hypothesis}
\numberwithin{equation}{section}

\newcommand{\eps}{\varepsilon}

\DeclareMathOperator{\id}{id} 
\DeclareMathOperator{\Ker}{Ker}
\DeclareMathOperator{\Ima}{Im}
\newcommand*\diff{\mathop{}\!\mathrm{d}}

\newcommand{\N}{\ensuremath{\mathbb{N}}}
\newcommand{\R}{\ensuremath{\mathbb{R}}}
\newcommand{\E}{\ensuremath{\mathcal{E}}}

\newcommand{\Kapp}{\ensuremath{\vec{\kappa}}}
\newcommand{\norm}[2][]{\ensuremath{\big\|#2\big\|_{#1}}}
\newcommand{\metric}[2][]{\ensuremath{\left\langle#2\right\rangle_{#1}}}
\newcommand{\betrag}[1]{\ensuremath{\left|#1\right|}}
\newcommand{\st}{\ensuremath{\,\big|\,}}
\newcommand{\deps}{\ensuremath{\partial_\eps}}
\newcommand{\diffs}{\ensuremath{\partial_s}}
\newcommand{\diffseps}{\ensuremath{\partial_{s_\eps}}}

\title{The \L ojasiewicz-Simon gradient inequality\\ for open elastic curves}
\author{Anna Dall'Acqua, Paola Pozzi, Adrian Spener}

\begin{document}

\maketitle

\begin{abstract}
In  this paper we consider the elastic energy for open curves in Euclidean space subject to clamped boundary conditions and obtain the \L ojasiewicz-Simon gradient inequality for this energy functional. Thanks to this inequality we can prove that a (suitably reparametrized) solution to the associated $L^2$-gradient flow converges for large time to an elastica, that is to a critical point of the functional.
\end{abstract}
 \bigskip
 


\section{Introduction}

In the past years a considerable number of papers dealing with the long-time existence of motion of curves  by the $L^{2}$-gradient flow for the  elastic energy have appeared in the literature. Closed  curves subject to a  length/area constraint of some sort and with or without  an inextensibility condition have been studied for instance in  \cite{DKS}, \cite{Polden}, \cite{Wen}, \cite{Koiso}, \cite{LangerSinger1}, \cite{Okabe2007}, \cite{Okabe2008};   open curves subject to different types of boundary conditions and a constraint on the length  have been treated in \cite{DP}, \cite{Lin}, \cite{DLP}, \cite{DLP2}, curves of infinite length in \cite{NO14}.

However, in several of the above frameworks the question of asymptotically convergence to an  equilibrium point of the gradient system has not been satisfactorily answered.
Indeed, in many of the works just mentioned (see for instance \cite{DKS}, \cite{DP}, \cite{Lin}, \cite{DLP}, \cite{DLP2}), the method of proof chosen to show long-time existence 
allows only to infer that for a sequence of time converging to infinity  there exists a subsequence of suitably reparametrized curves that converge to a critical point of the  energy functional. Thus, 
in principle different sequences could converge to different critical points.
Motivation for this work is to show  that this does not happen and that given an initial smooth curve the whole flow converges (after an appropriate reparametrization) to a stationary solution.

We give here a detailed proof for the setting  presented in \cite{Lin}, that is for the elastic flow of open curves subject to a constraint on the growth of the length (obtained by adding a suitable penalty term in the energy functional)  and clamped boundary conditions (i.e. the two boundary points of the curve and its tangents are kept fixed along the evolution).

More precisely, let us recall that
the elastic energy for a regular and sufficiently smooth curve $f:I \to \R^d$, $f=f(x)$, is given by
\begin{equation}\label{EE}
\E:f \mapsto \frac{1}{2}\int_I |\Kapp|^2 \diff s_f = \frac{1}{2}\int_I |\Kapp|^2 \,|f_x|\diff x,
\end{equation}
with $\Kapp$  the curvature vector, 
that is $\Kapp=\partial_{s_f} \partial_{s_f} f$ where $\partial_{s_f}= \frac{1}{|f_x|} \partial_x$. Here and in the following, $d \in \N$, $d\geq 2$, and $I :=[0,1]\subset \R$. It is well known that the energy $\E$ is a geometric functional, i.e. it is invariant under reparametrizations of the curve $f$, and that 
the $L^2$-gradient of the elastic energy is given by 
\begin{equation} \label{eq:FirstVariation}
\nabla_{L^2} \E (f) = \nabla_{s_f}^2 \Kapp + \frac{1}{2} |\Kapp|^2\Kapp,
\end{equation}
where $\nabla_{s_f} \phi : = \partial_{s_f} \phi - \langle \partial_{s_f} \phi ,\partial_{s_f} f \rangle \partial_{s_f} f$, see for instance \cite[Lemma A.1]{DP}.

Since the energy $\E$ might be decreased  by letting the curve grow towards infinity (just think of a (portion of a) circle whose radius is expanding), it is typical to penalize the growth of the length of the curve by considering the functional
\begin{align}\label{e-lambda}
 \mathcal{E}_\lambda(f) :=\E(f)+ \lambda \mathcal{L}(f)
 \end{align}
for a given positive constant $\lambda$. The $L^2$-gradient is then given by
\begin{equation}\label{gradlin} 
\nabla_{L^2} \E_\lambda (f) = \nabla_{s_f}^2 \Kapp + \frac{1}{2} |\Kapp|^2\Kapp -\lambda \Kapp,
\end{equation}
and the associated evolution reads
\begin{equation}
 \left \{ \begin{array}{ll}
           \partial_t f = -  \nabla_{s_{f}}^2 \Kapp - \frac{1}{2} |\Kapp|^2\Kapp +\lambda \Kapp & \;\; t \in (0,T),\\
           f(t,0) =f_-, \quad f(t,1)=f_+& \;\; t \in (0,T),\\
           \partial_{s_{f}} f (t,0) =T_-, \quad  \partial_{s_{f}} f (t,1) =T_+ & \;\; t \in (0,T),\\
           f(0,\cdot) = f_0(\cdot)&
          \end{array}
  \right.   \label{eq:ElasticFlow}
\end{equation}
for given $f_\pm \in \R^d$, unit vectors $T_\pm$, as well as a smooth regular initial curve $f_0:[0,1] \to \R^{d}$.  
The evolution problem \eqref{eq:ElasticFlow} has been studied in \cite{Lin}. There the following result is shown.
\begin{mythm}[{\cite[Theorem~1]{Lin}}]\label{ThmLin}
For any prescribed constant $\lambda \in (0, \infty)$ and smooth initial curve $f_0$ with finite length $\mathcal{L}(f_0) \in (0,\infty)$ and the clamped boundary conditions
\[f(t,0)=f_-, \quad f(t, 1)=f_+, \quad \partial_{s_{f}} f(t, 0)=T_-, \quad \partial_{s_{f}} f(t, 1)=T_+ \]
 there exists a global smooth solution to the $L^2$-flow of $\E_\lambda(f)$ in \eqref{e-lambda}. Moreover, after reparametrization by arc-length, the family of curves $f(t)$ subconverges to $f_\infty$, which is an equilibrium of the energy functional $\E_\lambda$.
\end{mythm}

In this work we strengthen the above statement by showing
\begin{mythm}\label{main2}
Let $\lambda > 0$ and $f:[0,\infty) \times I \to \R^d$ be the global smooth solution to the elastic flow \eqref{eq:ElasticFlow}. Then there exists a familiy of smooth diffeomorphisms $\Phi(t): I \to I$, $t \in (0,\infty)$ such that $f(t, \Phi (t,\cdot))$ converges  
to a critical point $\hat{f}_\infty$ of $\mathcal{E}_{\lambda}$, i.e.
\begin{equation}
 \norm[L^2]{f(t,\Phi (t,\cdot)) - \hat{f}_\infty} \to 0 \text{ as }t \to \infty. \label{eq:convergence}
\end{equation}
\end{mythm}

The reason for  focusing first on a framework dealing with  open curves  is that  a related statement has been proved in \cite{blowups} for \emph{closed surfaces} and we can expect to be able to adapt  those arguments to the case of closed curves without too much effort.  
More precisely, in \cite{blowups} the authors show that if the Willmore flow is started sufficiently close to a critical point, then the flow exists globally in time and it converges after a suitable reparametrization to a Willmore surface.

Our strategy in proving Theorem~\ref{main2} is to exploit Theorem~\ref{ThmLin} in order to identify a critical point of the energy functional and ``get sufficiently close'' to it and then employ some ideas from \cite{blowups} to show our claim.
In particular, inspired by \cite{blowups}, we also  apply  the  \L ojasiewicz-Simon gradient  inequality.  
In this respect a considerable effort in our work is spent in showing that the elastic energy functional satisfies the \L ojasiewicz-Simon gradient  inequality near a critical point, see Theorem~\ref{mainLos} below. The particular choice of boundary conditions imposes a particular choice of function spaces  and this in turn calls for new ideas and new arguments.

One of the advantages for using the \L ojasiewicz-Simon gradient inequality is that, under suitable circumstances, it provides  a $L^{1}$-control (in time) for  the velocity $\| \partial_{t} f \|_{L^{2}(I)} $ (see \eqref{derGest} and \eqref{decadenza} below), as opposed  to  the sole steepest descent property of the evolution, that yields only a $L^{2}$-control (see \eqref{eq:monotone}).
For  more motivation  and further applications of the \L ojasiewicz-Simon gradient inequality  we refer to \cite{Chill} and \cite{CF}.

Last but not least let us mention that in the case of \emph{planar} curves the statement of Theorem~\ref{main2} is obtained in \cite{NOpreprint} by different methods: there, starting from the long-time and sub-convergence result of Lin \cite{Lin},  the authors show that the set of critical points corresponding to any possible energy level is finite (up to reparametrization of the curve) and that these critical manifolds are isolated in the Hausdorff distance. Then, quite intuitively, the flow  has no other choice but to converge to one possible stationary solution.

The paper is organized as follows: after introducing important notation and discussing some preliminaries  in Section~\ref{sec2}, we devote Section~\ref{sec3} and \ref{sec4} to the proof of the \L ojasiewicz-Simon gradient inequality as given in Theorem~\ref{mainLos}. The proof of  Theorem~\ref{main2} is presented in Section~\ref{sec5}.
To improve the readability of the paper we have decided to collect lengthy and technical calculations   and auxiliary results in the Appendix.

\bigskip

\textbf{Acknowledgements:} Anna Dall'Acqua and Adrian Spener would like to thank Reiner Sch\"atzle for a very interesting and helpful discussion.
Paola Pozzi would like to thank the Isaac Newton Institute for Mathematical Sciences, Cambridge, for support and hospitality during the programme \emph{Coupling Geometric PDEs with Physics for Cell Morphology, Motility and Pattern Formation} where work on this paper was undertaken.


\section{Notation}\label{sec2}

First of all we provide some notation and recall some useful facts.
Euclidean scalar product and norm in $\R^d$ are denoted by $ \langle\cdot,\cdot \rangle_\text{euc}$ and $| \cdot|$ respectively. 
The specification \lq\lq$\text{euc}$\rq\rq ~will be omitted if clear from the context.
For simplicity of notation the dot \lq\lq$\cdot$\rq\rq ~will also refer to the Euclidean scalar product.
A scalar product in an arbitrary Hilbert space $H$ will be denoted by $\langle\cdot,\cdot \rangle_H$.
If no index is specified then $\langle\cdot,\cdot \rangle$ stands for the usual dual pairing with respect to a (specified) Banach space and its dual. A constant $C$ can always change from line to line.
 
In the following let $d \in \N$, $d\geq 2$, and $I =[0,1]\subset \R$. Moreover, let $f \in H^4(I,\R^d)$, $f=f(x)$, denote a regular curve. Regularity means that $|f_x| \neq 0$ in $I$ and hence that $f$ is an immersion. The curve $f$ induces the volume form $\diff s_f = |\partial_x f|\diff x$ on $I$. 
For vector fields $\phi, \psi: I \to \R^d$, $p \in [1, \infty)$, $\psi=\psi(x)$ and $\phi=\phi(x)$, we write
\begin{align*}
\|\psi \|_{L^p(\diff s_f)} := \Big(\int_I |\psi|^p  \diff s_f \Big)^{1/p} \quad \text{and} \qquad \langle \phi, \psi \rangle_{L^2 (\diff s_f) }= \int_I \phi \cdot \psi  \, \diff s_f.
\end{align*}
From the Sobolev embedding $H^2\hookrightarrow \mathcal{C}^1$ and the regularity of $f$  we know that $|\partial_x f| \in [\varepsilon,\frac{1}{\varepsilon}]$ uniformly on $I$ for some $\varepsilon > 0$. Thus the $L^p$-spaces with respect to $\diff x$ and $\diff s_f$ on $I$ coincide.

For given $\phi:I \to \R^{d}$ we denote by $\nabla_{s_f} \phi$  the normal component of $\partial_{s_f} \phi$, that is
\[\nabla_{s_f} \phi : = \partial_{s_f} \phi - \langle \partial_{s_f} \phi ,\partial_{s_f} f \rangle \partial_{s_f} f.\]
The (weak) derivative with respect to $s_f$ is simply defined
 by $\partial_{s_f} (\cdot)=\frac{1}{|f_x|} \partial_x (\cdot)$.
 Note that the definition is meaningful since for any $\eta:I \to \R^d$ and $\varphi \in C^\infty_0(I, \R^d)$ we have
 \begin{align*}
 \int_I  \eta \cdot \partial_{s_f} \varphi \diff s_f =\int_I \eta  \cdot \varphi_x \, \diff x = -\int_I \eta_x \cdot \varphi \diff x =
 -\int_I  \partial_{s_f} \eta \cdot \varphi \diff s_f.
 \end{align*}
 We also can define the Sobolev spaces $H^k(I,\diff s_f;\R^d)$ of mappings from $I$ to $\R^d$ which are squared-integrable, $k$-times weakly differentiable with weak \mbox{derivatives} in $L^2(I,\diff s_f; \R^d)$ and equipped with the norm $\norm[H^k (\diff s_f)]{\phi} : = \sum_{i=0}^k \norm[L^2(\diff s_f)]{\partial_{s_f}^i \phi}$. Note that for the definition of  $H^k(I,\diff s_f;\R^d)$ with $ k >4$ the regularity of $f$ must be increased. Moreover $H^k(I,\diff s_f;\R^d)$ and $H^{k}(I, \R^{d})$ are the same spaces: the specification of the metric and associated measure will be given only when necessary. We will work with both since when $f$ changes in the course of our arguments  it is convenient to stick to a fixed metric, that is to work with $H^{k}(I, \R^{d})$. The use of $H^k(I,\diff s_f;\R^d)$ will be preferred when we apply purely geometric arguments.

\begin{mydef} 
The space $H^{4}_{c}$ is defined as
\begin{equation}\label{eq:DefH^4Clampedall}
H^{4}_{c} := H^{4}(I, \R^{d}) \cap H^{2}_{0}(I, \R^{d}) \,. 
\end{equation}
\end{mydef}
Note that because of the above remarks we have that $H^4_c = H^4_c (\diff s_f)$, with $H^4_c (\diff s_f)= H^{4}(I, \diff s_f)\cap H^2_0(I,\diff s_f)$, for any regular curve $f \in H^4(I,\R^d)$.\bigskip 

Eventually we will work on the space of variations normal to a given fixed regular curve $\bar{f}$.  With this in mind we introduce the following spaces that depend on the particular choice of the given immersion $\bar{f}$.
\begin{mydef}\label{def2.2}
Let $\bar{f} \in H^{4}(I, \R^{d})$ be an immersion.
The subspace of $H^{4}_{c}$ given by normal variations is denoted by
\begin{equation}\label{eq:DefH^4Clamped}
H^{4,\bot}_{c} := \{\phi \in H^{4}_{c}  \st \metric[\text{euc}]{\phi, \partial_x \bar{f}} (x) = 0 \;  \forall x\in I \}. 
\end{equation}
Moreover
\begin{equation}
L^{2,\bot} := \{\phi \in L^{2}(I, \R^{d})  \st \metric[\text{euc}]{\phi, \partial_x \bar{f}} (x) = 0 \;  \forall x\in I\}. 
\end{equation}
\end{mydef}

In the \L ojasiewicz-Simon inequality we will describe the behavior of the elastic functional close to a critical point: first by studying only normal variations and then general ones. For this the following notation is useful.
\begin{mydef}\label{def2}
Let $\bar{f} \in H^{4}(I, \R^{d})$ be an immersion.
The restriction of the elastic energy \eqref{EE} to normal perturbations is defined as 
\begin{equation}\label{Esimple}
E:U  \to \R,\; \phi \mapsto \E(\bar{f}+\phi) \mbox{ with }U := B_\rho(0) \subset H^{4,\bot}_{c},
\end{equation}
while, for the ease of notation, we also define
\begin{equation}\label{Efat}
\mathbf{E}:\mathbf{U}  \to \R,\; \psi \mapsto \E(\bar{f}+\psi) \mbox{ with }\mathbf{U} := B_\rho(0) \subset H^{4}_{c}.
\end{equation}
Here $\rho>0$ is chosen small enough such that $\bar{f}+\psi$, respectively $\bar{f}+\phi$, is still immersed for all choices of $\psi \in \mathbf{U}$, respectively $\phi \in U$.
\end{mydef}

The functionals $E$ and $\mathbf{E}$ are Fr\'echet-differentiable with derivatives
$E\rq{}: U \to (H^{4,\bot}_{c})^{*}$ and $\mathbf{E}\rq{}: \mathbf{U} \to (H^{4}_c)^{*}$.
It turns out that these Fr\'echet-derivatives may be identified with (geometrical) $L^{2}$-gradients in the following way.
 The operator $\nabla E: U \to L^{2,\bot}$ is defined by the identity
\begin{equation*}
 E\rq{}(\phi) \eta = \langle \nabla E(\phi), \eta \rangle _{L^2(\diff s_{(\bar{f}+\phi)})} \; \forall \phi \in U, \forall \eta \in H^{4,\bot}_{c} \, .
\end{equation*}
Similarly, $\nabla \mathbf{E}: \mathbf{U} \to L^{2}$, $\mathbf{U} =B_\rho(0)\subset H^{4}_{c}$, is  defined by
\begin{equation*}
 \mathbf{E}\rq{}(\psi)\eta = \langle \nabla \mathbf{E}(\psi), \eta \rangle _{L^2 (\diff s_{(\bar{f} +\psi)})} \; \forall \psi \in \mathbf{U}, \forall \eta \in H^{4}_{c} \, .
\end{equation*}
Standard computations (see for instance \cite[Lemma A.1]{DP}) together with the fact that $\partial_x \eta =0$ at the boundary give that
\begin{align}\label{gradfat}
\nabla \mathbf{E}(\psi)= \nabla_s^2 \Kapp + \frac{1}{2} |\Kapp|^2 \Kapp,
\end{align}
where here $\Kapp$ is the curvature vector of the curve $\bar{f}+\psi$ and $s=s_{(\bar{f}+\psi)}$.
For $\nabla E(\phi)$ we get
\begin{align*}
\nabla E(\phi)&=\nabla \mathbf{E}(\phi)-  \frac{\bar{f}_{x}}{|\bar{f}_{x}|} \langle \nabla \mathbf{E}(\phi), \frac{\bar{f}_{x}}{|\bar{f}_{x}|} \rangle_\text{euc} \\
& = \nabla_s^2 \Kapp + \frac{1}{2} |\Kapp|^2 \Kapp -  \frac{\bar{f}_{x}}{|\bar{f}_{x}|} \langle \nabla_s^2 \Kapp + \frac{1}{2} |\Kapp|^2 \Kapp, \frac{\bar{f}_{x}}{|\bar{f}_{x}|} \rangle_\text{euc} ,
\end{align*}
where here $\Kapp$ is the curvature vector of the curve $\bar{f}+\phi$ and $s=s_{(\bar{f}+\phi)}$. The difference in the two expression is due to the fact that $\nabla E(\phi) \in L^{2,\bot}$, whereas $\nabla \mathbf{E}(\psi) \in L^2$.

For more clarity in the application of some functional analytic arguments  employed below, it is convenient to  identify  the Fr\'echet-derivative with the gradient with respect to a fixed $L^2$-metric. 
 More precisely we write
\begin{equation*}
 \mathbf{E}\rq{}(\psi)\eta = \langle d \mathbf{E}(\psi), \eta \rangle _{L^2 (\diff x)} \; \forall \psi \in \mathbf{U}, \forall \eta \in H^{4}_{c} \, .
\end{equation*}
The above observations give that  
\[
d \mathbf{E}(\psi)=\nabla \mathbf{E}(\psi) |\bar{f}_x + \psi_x|= (\nabla_s^2 \Kapp + \frac{1}{2} |\Kapp|^2 \Kapp) |\bar{f}_x + \psi_x|  
\]
with $s=s_{(\bar{f}+\psi)}$ and $\Kapp$ the curvature of $\bar{f}+\psi$.
The $L^2$-gradient $dE$ for normal variations is defined analogously, precisely
\[ dE(\phi)= \nabla E(\phi) |\bar{f}_x + \phi_x|=  \Big (\nabla_s^2 \Kapp + \frac{1}{2} |\Kapp|^2 \Kapp -  \frac{\bar{f}_{x}}{|\bar{f}_{x}|} \langle \nabla_s^2 \Kapp + \frac{1}{2} |\Kapp|^2 \Kapp, \frac{\bar{f}_{x}}{|\bar{f}_{x}|} \rangle_\text{euc}  \Big ) |\bar{f}_x + \phi_x|, \]
where $\Kapp$ is the curvature vector of the curve $\bar{f}+\phi$ and $s=s_{(\bar{f}+\phi)}$.

We are now ready to state the \L ojasiewicz-Simon inequality that we intend to prove.
\begin{mythm}\label{mainLos}
Let $\bar{f} \in H^{4}(I, \R^{d})$ be a  regular critical point for $\E$, i.e. $\nabla_{L^2}\E(\bar{f})=0$.
Then there exists constants $C_1 >0$, $\theta \in (0,\frac{1}{2}]$ and $\sigma > 0$ such that all $\psi \in B_\sigma(0) \subset \mathbf{U} \subset H^{4}_{c}$ satisfy the inequality
\begin{equation}\label{eq:LS}
 \betrag{\mathbf{E}(0)-\mathbf{E}(\psi)}^{1-\theta}   
 \leq C_1 \norm[L^2(\diff s_{(\bar{f}+\psi)})]{\nabla \mathbf{E} (\psi)} .
\end{equation}
\end{mythm}
First of all notice that 
since $\bar{f} \in H^4(I; \R^d)$ is a regular curve and a critical point of the elastic energy, it is smooth, as  one can see  for instance by a bootstrap argument (recall \eqref{eq:FirstVariation} and \eqref{eqa4}). In particular, $\bar{f} \in H^5(I; \R^d)$, which is the highest regularity assumption that will be employed in the arguments that follow.

We prove \eqref{eq:LS} in two steps. First, we show the inequality only for normal variations, that is for the functional $E$ (see Theorem~\ref{LosNormal} below. The main ideas for the proof are explained in Remark \ref{rem:list}). Then, we deduce from this the inequality in its general form. 


\section{Proof of the \L ojasiewicz-Simon inequality for normal variations}
\label{sec3}

The main result of this section is the following.
\begin{mythm}\label{LosNormal}
Let $\bar{f} $  be a smooth regular critical point for $\E$, i.e. $\nabla_{L^{2}}\E(\bar{f})=0$. Then there exists constants 
$C_2> 0$, $\theta \in (0,\frac{1}{2}]$ and $\tilde{\sigma} > 0$ such that all $\phi \in B_{\tilde{\sigma}}(0) \subset U \subset H^{4,\bot}_{c}$ satisfy the inequality
\begin{equation}\label{eq:LSNormal}
 \betrag{E(0)-E(\phi)}^{1-\theta} \leq C_2 \norm[L^2(dx)]{d E (\phi)}  .
\end{equation}
\end{mythm}

\begin{remark}\label{remgeometric}
Notice that in \eqref{eq:LS} only geometric quantities appear, whereas in Theorem \ref{LosNormal} we work with a fixed parametrization obtaining an estimate with $d E (\phi)$.
This reflects the fact that the proof of the \L ojasiewicz-Simon inequality for normal variations is mainly based on functional analytic arguments, whereas Thereom~\ref{mainLos} will be applied in a geometric context.
However, observe that
by Theorem \ref{LosNormal}, using that $\|\bar{f}_x+ \phi_x\|_{L^\infty} \leq C(\bar{f}, \tilde{\sigma})$ and the explicit expressions for $d E(\phi)$ and $\nabla E (\phi)$ given below Definition~\ref{def2}, we can find a constant $C_{3}>0$ such that with the same $\tilde{\sigma}$ and  $\theta \in (0,\frac{1}{2}]$ we have
\begin{align*}
\betrag{E(0)-E(\phi)}^{1-\theta} \leq C_2 \norm[L^2(dx)]{d E(\phi)}\leq C_{3} \| \nabla E (\phi) \|_{L^{2}(\diff s_{\bar{f}+ \phi})} 
\end{align*}
for all $\phi \in B_{\tilde{\sigma}}(0) \subset U \subset H^{4,\bot}_{c}$.
\end{remark}

\begin{remark}\label{rem:list} 
The  \L ojasiewicz-Simon inequality for normal variations follows once we have proven that for some $\rho >0$ small enough on $U=B_{\rho}(0) \subset  H^{4,\bot}_{c}$:
\begin{enumerate}
 \item $E:U\to \R$, $\phi \mapsto \E(\bar{f}+\phi)$ is analytic,
 \item its gradient $d E:U \to L^{2,\bot}$ is analytic and 
 \item the Fr\'echet derivative of the gradient $\mathfrak{L}:=(d E)'(0): H^{4,\bot}_{c} \to L^{2, \bot}$ evaluated at $0$ is a Fredholm operator. 
\end{enumerate}
This follows from \cite[Corollary 3.11]{Chill}. (See also \cite[Page 355]{blowups}.)  For completeness in Appendix \ref{Chill} we give the statement of \cite[Corollary 3.11]{Chill} and show that \textit{1., 2., 3.} above are sufficient  
for the assumptions in  \cite[Corollary 3.11]{Chill} to be satisfied. 
Note that the argument of Chill is more general since it  
relies merely on the structure of a Banach space not necessarily endowed with an inner product.
\end{remark}

\begin{proof}[Proof of Theorem \ref{LosNormal}]
By Remark \ref{rem:list}, the claim in Theorem \ref{LosNormal} follows from Theorem \ref{thmanaly} and Corollary \ref{cor:Fredholm} below.
\end{proof}


\subsection{Analyticity}

The aim of this section is to show the analyticity of $E$ and its gradient $d E$ (cf.~ Remark~\ref{rem:list}). For convenience we recall the definition and fundamental properties of analytic maps between Banach spaces, see \cite[Definition 8.8]{Zeidler} or \cite{analyticBanach}. Let $(X, \| \cdot \|_{X})$,  $(Y, \| \cdot \|_{Y})$, $(Y_1, \| \cdot \|_{Y_1})$, $(Y_2, \| \cdot \|_{Y_2})$,  $(Z, \| \cdot \|_{Z})$ be arbitrary  Banach spaces.
A map $f: D \to Y$, $D$ open in $X$, is \textit{analytic at $x_0 \in D$} if in a  neighborhood of $x_0$
\begin{align*}
 \sum_{k=0}^{\infty} \| a_k\| \|x-x_0\|^k \mbox{ converges} 
\mbox{ and } f(x)= \sum_{k=0}^{\infty} a_k (x-x_0)^k\, .
\end{align*}
Here each $a_k$ is a $k$-linear, symmetric and continuous map from $X^k= X \times \cdots \times X$ ($k$-times) into $Y$ and writing $a_k(x-x_0)^k$ we actually mean $a_k(x-x_0, ..., x-x_0)$. The map $f$ is \textit{analytic in $D$} if it is analytic at each point of $D$.

A linear and bounded map $F \in \mathcal{B}(X, Y)$ is analytic  and so is its affine counterpart $F(x_{1}+ \cdot)$ for a given $x_{1} \in X$. This follows immediately, since we may write $F(x)= F(x_{0})+ F(x-x_{0})= \sum_{k=0}^{\infty} a_{k} (x-x_{0})^k  $, where $a_{k}=0$ for $k>2$, $a_{1}=F \in \mathcal{B}(X, Y)$ and $a_0=F(x_0)$. With similar arguments one immediately sees that a map $F: X \to Y$, $x \to a_{2}(x,x)$ with $a_{2}: X \times X \to Y$ a bilinear, symmetric, continuous map, is analytic.  
Compositions of analytic functions is again analytic, that is if $G:D \to E$ and $F : E \to Z$ are analytic, where $D $ open in $X$, $E$ open in $Y$, then $F \circ G: D \to Z$ is analytic (see for instance \cite[page 1079]{analyticBanach}). If the Banach space $Y$ is continuosly embedded into a Banach space $Z$, then any analytic function from $D$, open set in $X$, to $Y$ is also analytic as a map from $D$ to $Z$. Moreover, the sum of analytic functions is clearly analytic.

In the following, it will be convenient to characterise the analyticity of a mapping seeing it as a product of two mappings. Since we work in Banach spaces we need some additional structure. If $F: D \to Y_1$, $G: D \to Y_2$, $D $ open in $X$, are analytic and there exists a Banach space $Z$ as well as a bilinear continuous mapping $*:Y_1 \times Y_2 \to Z$ then the product $F*G:D \to Z$, $x \mapsto F(x)*G(x)$, is analytic. This can be proved using similar ideas as for the Cauchy product of series. Consequently the same holds for products with a finite number of factors. A particular case is when  
$Y_1=Y_2=Y$ and $Y$ is a Banach algebra. In this case, the product of two analytic functions $G,F: D \to Y$  
 is again analytic as a mapping from $D$ to $Y$. In the proof below we will use that $H^m(I,\R)$, $m\geq 1$, is a Banach algebra (see \cite[Cor.8.10]{Brezis} for $m=1$).

\begin{lemma}\label{lemma:analycitybyparts}
Let $\bar{f} \in H^{4}(I;\R^d)$ be a regular curve and $U=B_{\rho}(0) \subset H^{4,\bot}_{c}$ be as in Definition \ref{def2}.
The following functions are well defined and (real-) analytic:
\begin{enumerate}
 \item  $F_1 :U \to  H^3(I,\R)$, $\phi \mapsto |\partial_x (\bar{f} + \phi)|$,
  \item $F_2: U \to H^3(I, \R^{d}), \phi \mapsto \frac{\partial_x (\bar{f} + \phi)}{|\partial_x (\bar{f} + \phi)|}= \partial_s (\bar{f} + \phi)$ where $s =s_{(\bar{f}+ \phi)}$,
\item $F_3 : U \to H^2(I, \R^{d}), \phi \mapsto \Kapp_{\bar{f}+\phi}= \partial_s^2 (\bar{f} + \phi)$ where $s =s_{(\bar{f}+ \phi)}$,
  \item $F_4 : U \to L^1(I, \R)$, $ \phi \mapsto |\Kapp_{\bar{f}+\phi}|^2 |\partial_x (\bar{f} + \phi)|$,
\item $F_5: U \to L^{2}$, 
  $\phi \mapsto   \Big (\nabla_s^2 \Kapp + \frac{1}{2} |\Kapp|^2 \Kapp  \Big ) |\bar{f}_x + \phi_x|$ , where $\Kapp= \Kapp_{\bar{f}+\phi}$, $s =s_{(\bar{f}+ \phi)},$
\end{enumerate}
\end{lemma}

For sake of readability the proof of Lemma \ref{lemma:analycitybyparts} is given in Appendix \ref{appanaly}. 
\begin{mythm}\label{thmanaly}
Let $\bar{f} \in H^{4}(I;\R^d)$ be a regular curve and $U=B_{\rho}(0) \subset H^{4,\bot}_{c}$ be as in Definition \ref{def2}.
Then the functions 
\[E : U \to \R, \;  \phi \mapsto \mathcal{E}(\bar{f}+\phi) =\frac{1}{2}\int_I |\Kapp_{\bar{f}+\phi}|^2 \diff s_{\bar{f}+\phi}\]
and 
\[d E: U \to L^{2,\bot}, \; \phi \mapsto   \Big (\nabla_s^2 \Kapp + \frac{1}{2} |\Kapp|^2 \Kapp -  \frac{\bar{f}_{x}}{|\bar{f}_{x}|} \langle \nabla_s^2 \Kapp + \frac{1}{2} |\Kapp|^2 \Kapp, \frac{\bar{f}_{x}}{|\bar{f}_{x}|} \rangle_\text{euc}  \Big ) |\bar{f}_x + \phi_x| \, ,\]
(where $\Kapp= \Kapp_{\bar{f}+\phi}$, $s =s_{(\bar{f}+ \phi)}$) are (real-) analytic.
\end{mythm}
\begin{proof}
The analyticity of the function $E$ follows directly from the analyticity of the function $F_4$ defined in Lemma  \ref{lemma:analycitybyparts} since integration over $[0,1]$ is a well-defined continuous and linear operator on $L^1(I,\R)$. Since the projection 
\[F: L^2(I,\R^d) \to  L^{2,\bot}(I,\R^d), \quad \phi \mapsto \phi - \frac{\bar{f}_{x}}{|\bar{f}_{x}|} \langle \phi, \frac{\bar{f}_{x}}{|\bar{f}_{x}|} \rangle_\text{euc}\, , \]
is a linear and continuous operator,  the analyticity of $dE$ is a direct consequence of the analyticity of the function $F_5$ defined in Lemma  \ref{lemma:analycitybyparts}.
\end{proof}


\subsection{The second variation of \texorpdfstring{$E$}{E}}

\begin{prop}\label{prop:2ndvar}
Let $\bar{f} \in H^{4}(I;\R^d)$ be a regular curve and $\phi, \psi \in U \subset H^{4,\bot}_\text{c}$. Then the second variation of $E$ (defined in \eqref{Esimple}) at $0$ in the direction of $\phi$ and $\psi$ is given by
 \begin{align} \label{eq:2ndvar}
 \begin{split} E''(0)(\psi,\phi)  &= \int_I \left(\metric{\nabla_s^2 \phi, \nabla_s^2 \psi} +\metric{ \nabla_s^2 \phi,\Kapp}\metric{\psi,\Kapp} +\metric{\nabla_s^2 \psi,\Kapp} \metric{\phi,\Kapp}\right. \\
&\qquad  -\frac{3}{2}|\Kapp|^2\metric{\nabla_s \psi, \nabla_s \phi}+ \metric{\nabla_s \psi, \Kapp}\metric{\nabla_s \phi,\Kapp}\\
&\qquad   - \metric{\nabla_s\psi,\nabla_s\Kapp} \metric{\phi, \Kapp}  -\metric{\nabla_s\phi,\nabla_s\Kapp} \metric{\psi, \Kapp}\\
&\qquad+ \left. |\Kapp|^2 \metric{\Kapp, \psi} \metric{\Kapp,\phi}\right) \diff s_{\bar{f}}     ,      \end{split}
\end{align}
where
$\metric{\cdot,\cdot }= \langle\cdot,\cdot\rangle_{\text{euc}}$, $s=s_{\bar f}$ and $\Kapp=\Kapp_{\bar{f}}$.
\end{prop}

The lengthy calculation of the second variation is given in \ref{proof:2ndvar} in the appendix.

We also immediately find the second variation of the elastic energy with penalized length $E_\lambda:U\to \R$, $\phi \mapsto E(\phi) + \lambda \mathcal{L}(\bar{f}+\phi)$ for $\lambda \in \R, \lambda \geq 0$.
\begin{corol}\label{cor:2varl}
Under the same assumptions as in Proposition \ref{prop:2ndvar}, the second variation of $E_\lambda$ at $0$ in the direction of $\phi$ and $\psi$ is given by
 \begin{align} \label{eq:2ndvarlambda}
 \begin{split} E_\lambda''(0)[\psi,\phi]  &= \int_I \left( \metric{\nabla_s^2 \phi, \nabla_s^2 \psi} +\metric{ \nabla_s^2 \phi,\Kapp}\metric{\psi,\Kapp} +\metric{\nabla_s^2 \psi,\Kapp} \metric{\phi,\Kapp}\right. \\
&\qquad  +(\lambda-\frac{3}{2} |\Kapp|^2)\metric{\nabla_s \psi, \nabla_s \phi}+ \metric{\nabla_s \psi, \Kapp}\metric{\nabla_s \phi,\Kapp}\\
&\qquad   - \metric{\nabla_s\psi,\nabla_s\Kapp} \metric{\phi, \Kapp}  -\metric{\nabla_s\phi,\nabla_s\Kapp} \metric{\psi, \Kapp}\\
&\qquad+ \left. |\Kapp|^2 \metric{\Kapp, \psi} \metric{\Kapp,\phi}\right) \diff s_{\bar f},          \end{split}
\end{align}
where
$\metric{\cdot,\cdot }= \langle\cdot,\cdot\rangle_{\text{euc}}$, $s=s_{\bar f}$ and $\Kapp=\Kapp_{\bar{f}}$.
\end{corol}

 This calculation is given in \ref{proof:2ndvarlambda} in the appendix.


\subsection{The Fredholm property}
In this section we show that the operator $\mathfrak{L}$ associated to the second variation of $E$ (defined in \eqref{Esimple}) in $0$ is Fredholm of index $0$ (cf.~ Remark~\ref{rem:list}). To do so, we derive the Fredholm property for  the leading term (which is associated to the  bilinear form defined in Definition \ref{def:bili} below) and then use compact embedding theorems for the perturbation. More precisely, we first show the following result:
\begin{prop}\label{prop:FredholmMain}
Let $\bar{f} \in H^{5}(I;\R^d)$ be a regular curve. Then the operator
\begin{equation}
\nabla_s^4: H^{4,\bot}_{c} (\diff s_{\bar f})
\to L^{2,\bot} (\diff s_{\bar f})
\end{equation}
is Fredholm of index $0$.
\end{prop}

We start by considering the following bilinear form.
\begin{mydef}\label{def:bili}
 Let $\bar{f}$ be as in Definition~\ref{def2.2} and denote the subspace of normal vector fields by
\[H^{2, \bot}_{0}(\diff s_{\bar f})= H^{2}_{0}((I,\diff s_{\bar f}), \R^{d})\cap L^{2, \bot}(\diff s_{\bar f})\, .\]
We let $a_{\bar{f}}$ be the form given by
 \[
  a_{\bar{f}}:H^{2, \bot}_0(\diff s_{\bar f})\times H^{2, \bot}_0 (\diff s_{\bar f})\to\R,\;\; (\phi,\psi) \mapsto \int_I \nabla^{2}_s \phi \cdot \nabla^{2}_s \psi \; \diff s_{\bar{f}}.
  \]
 Here $s=s_{\bar{f}}$.
\end{mydef}

\begin{lemma} \label{lemma:ellipticity}
 The bilinear form $a_{\bar{f}}$ defined above is bounded, symmetric and $H^{2,\bot}_0$-elliptic, i.e. there exist constants $\omega,\mu > 0$ such that 
 \begin{equation}
  a_{\bar{f}}(\phi,\phi) + \omega \norm[L^{2}(\diff s_{\bar{f}})]{\phi}^2 \geq \mu\norm[H^{2}(\diff s_{\bar{f}})]{\phi}^2 \text{ for all }\phi \in H^{2, \bot}_{0}(\diff s_{\bar f}).
 \end{equation}
\end{lemma}

\begin{proof}  Symmetry and boundedness are straight forward.
Moreover notice that, due to 
\begin{align}\label{uu1}
\int_{I} |\partial_{s} \phi|^{2} \diff s_{\bar{f}}= -\int_{I} \langle \phi, \partial_{s}^{2} \phi \rangle_{\text{euc}}\diff s_{\bar{f}}
\end{align} and the Cauchy-Schwarz inequality, the norm $\|\phi\|_{H^{2}(\diff s_{\bar{f}})}$ is equivalent to the norm $\| \phi \|:=\|\phi \|_{L^{2}(\diff s_{\bar{f}})} + \|\partial_{s}^{2} \phi\|_{L^{2}(\diff s_{\bar{f}})}$.
 Here and in the following $s=s_{\bar{f}}$ and $\Kapp=\Kapp_{\bar{f}}$ (and we omit the index $\bar{f}$ for the sake of simplicity of notation). 
Using that $\|\Kapp\|_{L^{\infty}}, \|\partial_s \Kapp\|_{L^{\infty}} \leq C(\bar{f})$, \eqref{eqa4}, \eqref{uu1} and Young's inequality it follows that for any $\varepsilon>0$
\begin{align*}
 a_{\bar{f}}(\phi,\phi) = \int_{I}  |\nabla_{s}^{2} \phi |^{2} \diff s_{\bar f} & \geq  \int_{I}  (\frac12 |\partial_{s}^{2} \phi|^{2} - C | \partial_{s} \phi |^{2} - C  |\phi|^{2}) \diff s_{\bar f}\\
& \geq \int_{I} \frac12 |\partial_{s}^{2} \phi |^{2} \diff s_{\bar f} - \varepsilon C  \int_{I} |\partial_{s}^{2} \phi |^{2} \diff s_{\bar f} - C_{\varepsilon}  \int_{I} |\phi|^2 \diff s_{\bar f} \, ,
\end{align*}
with a constant $C_{\varepsilon}$ depending only on $\bar{f}$ and $\varepsilon$.
Hence, for any $\omega > 0$
\begin{align*}
 a_{\bar{f}}(\phi,\phi) + \omega \norm[L^2(\diff s_{\bar{f}})]{\phi}^2 
 &\geq \int_I (\frac12 - C\varepsilon) |\partial_{s}^{2} \phi|^2 + \left(\omega - C_{\varepsilon} \right)|\phi|^2 \diff s_{\bar f}. 
 \end{align*}
Choosing $\varepsilon$ small enough and subsequently $\omega$ large we find that the last term is larger than $\frac{1}{4} (\| \phi \|^{2}_{L^{2}(\diff s_{\bar{f}})} + \| \partial_{s}^{2 }\phi \|^{2}_{L^{2}(\diff s_{\bar{f}})} )$, which proves the claimed ellipticity.
\end{proof}

\begin{lemma}\label{lem2.12}
Let $a_{\bar{f}}$ be the bilinear form defined in Definition \ref{def:bili} and $\omega$ as in Lemma \ref{lemma:ellipticity}. Then there exists a unique operator $A: D(A) \to L^{2,\bot}(\diff s_{\bar{f}})$ such that
 \begin{align} \label{eq:DefA}
D(A) &= \{ \phi \in H^{2,\bot}_0 (\diff s_{\bar{f}})\st \exists\, \xi \in L^{2,\bot} \text{ s.t. }\\ \nonumber
& \qquad  a_{\bar{f}}(\phi, \psi)  =  \langle \xi, \psi \rangle_{L^{2,\bot}(\diff s_{\bar{f}})} \text{ for all } \psi \in H^{2,\bot}_0 (\diff s_{\bar{f}}) \}\\ \label{eq:DefAA}
A\phi &= \xi.
\end{align}
The domain $D(A)$ with norm $\| \phi \|_{D(A)}:= \| \phi \|_{H^{2}(\diff s_{\bar{f}})} + \| A \phi \|_{L^2(\diff s_{\bar{f}})}$ is a Banach space and
the operator $A+\omega I$ is an isomorphism from $D(A)$ to $L^{2,\bot}(\diff s_{\bar{f}})$. \end{lemma}
\begin{proof} 
By Lemma \ref{lemma:ellipticity} the bilinear form $a_{\bar{f}}(\cdot, \cdot) + \omega \langle \cdot, \cdot \rangle_{L^2(ds_{\bar{f}})}$ on $H^{2,\bot}_0(\diff s_{\bar{f}})$ is bounded, symmetric and $H^{2,\bot}_{0}$-elliptic. Hence, for any $h \in L^{2,\perp}(\diff s_{\bar{f}})$ by Lax-Milgram Theorem there exists a unique $\phi \in H^{2,\bot}_0(\diff s_{\bar{f}}) $ such that
\begin{align}\label{eqim}
a_{\bar{f}}(\phi, \psi) + \omega \langle \phi, \psi \rangle_{L^2(ds_{\bar{f}})} =  \langle h, \psi \rangle_{L^2(ds_{\bar{f}})}  \; \mbox{ for all }\psi \in H^{2,\bot}_0 (\diff s_{\bar{f}}).
\end{align}
This defines an injective and linear operator 
\[\tilde{A}_{\omega}: L^{2,\bot} \to H^{2,\bot}_0,  \quad \tilde{A}_{\omega} (h) = \phi \mbox{ with }\phi\mbox{ solving \eqref{eqim}.}\]
The continuity of the operator follows from its coercivity and \eqref{eqim} since
\begin{align*}
 \|\tilde{A}_{\omega}(h) \|^{2}_{H^{2}(\diff s_{\bar{f}})}  & =\| \phi \|^{2}_{H^{2}(\diff s_{\bar{f}})} \leq \frac{1}{\mu} ( a_{\bar{f}}(\phi, \phi) + \omega \langle \phi, \phi \rangle_{L^2(ds_{\bar{f}})} ) \\
& = \frac{1}{\mu} \langle h, \phi \rangle_{L^2(ds_{\bar{f}})} \leq \frac{1}{\mu}\| h \|_{L^2(ds_{\bar{f}})} \| \phi \|_{H^2(ds_{\bar{f}})} \\
&= \frac{1}{\mu}\| h \|_{L^2(ds_{\bar{f}})} \| \tilde{A}_{\omega}(h) \|_{H^2(ds_{\bar{f}})} , 
\end{align*}
with $\mu$ as in Lemma \ref{lemma:ellipticity}. 

Then $\tilde{A}_{\omega}$ is a linear and continuous bijection on its range $R(\tilde{A}_{\omega}) \subset  H^{2,\bot}_0(\diff s_{\bar{f}})$, where
\begin{align*}
R(\tilde{A}_{\omega}) & =  \{\phi \in H^{2,\bot}_0 \st \exists\, \xi \in L^{2,\bot} \text{ s.t. } \\
& \qquad \quad \qquad \quad a_{\bar{f}}(\phi, \psi)  + \omega \langle \phi, \psi \rangle_{L^2(ds_{\bar{f}})}  =  \langle \xi, \psi \rangle_{L^{2}(\diff s_{\bar{f}})} \; \forall  \psi \in H^{2,\bot}_0 (\diff s_{\bar{f}})\} \, .
\end{align*} 
 Let 
\[A_{\omega}:= (\left.  \tilde{A}_{\omega}\right|_{R(\tilde{A}_{\omega})})^{-1}: D(A_{\omega}) \rightarrow L^{2,\perp}(\diff s_{\bar{f}})\]
with $D(A_{\omega}) = R(\tilde{A}_{\omega}) \subset  H^{2,\bot}_0$ and such that for $\phi \in  D(A_{\omega})$
\[ \langle A_{\omega} (\phi), \psi \rangle_{L^2(ds_{\bar{f}})}  = a_{\bar{f}}(\phi, \psi) + \omega \langle \phi, \psi \rangle_{L^2(ds_{\bar{f}})}  \; \mbox{ for all }\psi \in H^{2,\bot}_0 (\diff s_{\bar{f}}).\]
One immediately sees that $A_{\omega}$ is a closed operator. Indeed, let $(x_n)_{n \in \N}$  be a sequence in $D(A_{\omega})$ with $x_n \to x$ in $ H^{2,\bot}_0$ and $A_{\omega} x_n \to y $ in $L^{2,\perp}$. Let $y_n:= A_{\omega} x_n \in  L^{2,\perp}$ for all $n \in \N$. By definition of $A_{\omega}$, $\tilde{A}_{\omega} y_n = x_n$. Since $y_n \to y $ in $L^2$ and $\tilde{A}_{\omega}$ is continuous, it follows that $\tilde{A}_{\omega}y_n \to \tilde{A}_{\omega}y =: \tilde{x} $ in $H^{2,\bot}_0$. The closedness of $A_{\omega}$ follows directly once we have shown that $x=\tilde{x}$ in $H^{2,\bot}_0$. This follows since for any $\varepsilon>0$ we have
\[ \| x - \tilde{x}\|_{H^2(ds_{\bar{f}})} \leq \| x - x_n\|_{H^2(ds_{\bar{f}})} + \| x_n - \tilde{x}\|_{H^2(ds_{\bar{f}})} < \frac12 \varepsilon + \| \tilde{A}_{\omega} y_n - \tilde{A}_{\omega}y\|_{H^2(ds_{\bar{f}})} < \varepsilon \, , \]
for $n$ sufficiently big.

Since $A_{\omega}$ is a closed operator, $D(A_{\omega})$ with norm $\| \phi \|_{H^{2}} + \| A_{\omega} \phi \|_{L^2}$ is a Banach space and the operator $A_{\omega}$ is an isomorphism from $D(A_{\omega})$ 
to $L^{2,\perp}$, by the open mapping theorem. We observe that
\begin{align*}
D(A_{\omega})= R(\tilde{A}_{\omega}) = D(A) \mbox{ as defined in \eqref{eq:DefA}}\, .
\end{align*} 
The claim follows considering the operator $A: D(A) \to L^{2,\perp}$ acting as follows $A\phi = A_{\omega} \phi - \omega \phi $ for all $\phi \in D(A)$ and observing that on $D(A)$ the two norms
\[ \phi \mapsto  \| \phi \|_{H^{2}(ds_{\bar{f}})} + \| A_{\omega} \phi \|_{L^2(ds_{\bar{f}})} \mbox{ and }
 \phi \mapsto  \| \phi \|_{H^{2}(ds_{\bar{f}})} + \| A \phi \|_{L^2(ds_{\bar{f}})}\]
are equivalent.
\end{proof}


We now characterise the domain of the operator $A$  defined in the Lemma \ref{lem2.12}. Precisely, we show that $D(A)=H^{4, \bot}_c(\diff s_{\bar{f}})$ (see \eqref{eq:DefH^4Clamped}).

\begin{lemma}\label{lem2.13}
Let $A$ be the operator defined in the Lemma \ref{lem2.12} and assume that $\bar{f} \in H^{5}(I;\R^d)$ is a regular curve. Then, $ D(A)=H^{4, \bot}_c(\diff s_{\bar{f}})$, $\| \cdot \|_{D(A)}$ and $\| \cdot \|_{H^4(\diff s_{\bar{f}})}$ are equivalent norms and (with $s=s_{\bar f}$) $A \phi= \nabla_s^4 \phi$ for all $\phi \in D(A)$.
\end{lemma}
\begin{proof}
Observe that if $\phi \in H^{4, \bot}_{c}$, then $ \phi \in H^{2,\bot}_0$ and   
\[ a_{\bar{f}}(\phi, \psi)= \int_{I} \nabla_{s}^{4} \phi \cdot \psi \, \diff s_{\bar{f}} = \langle \nabla_{s}^{4} \phi, \psi \rangle_{L^{2}(\diff s_{\bar{f}})} \mbox{ for all }\psi \in H^{2,\bot}_0 (\diff s_{\bar{f}}), \]
with $\nabla_{s}^{4} \phi \in L^{2, \bot}$. Then $\phi \in D(A)$ and 
\begin{equation}\label{action}
A \phi = \nabla_s^4 \phi \mbox{ for all }\phi \in H^{4, \bot}_{c}(\diff s_{\bar{f}}) \, .
\end{equation}

We prove now the other inclusion, namely $D(A) \subset H^{4, \bot}_{c}$. Since $D(A) \subset H^{2,\bot}_0$, we only need to show that any $\phi \in D(A)$ admits weak derivatives of order three and four and that these are in $L^2$. It is convenient here for a vector field $\psi: I \to \R^d$ to write $\psi = \psi^{\top}+ \psi^{\perp}$ with $\psi^{\top}:= \langle \psi, \partial_s \bar{f} \rangle_\text{euc} \partial_s \bar{f}$. We show first that for $\phi \in D(A) $
\begin{align}\label{ziel1}
 \left| \int_{I} \partial_{s}^{2} \phi \cdot \partial_{s} \psi \diff s_{\bar{f}} \right|  \leq C \| \psi\|_{L^2(\diff s_{\bar{f}})} \mbox{ for all }\psi \in H^{1}_0 (\diff s_{\bar{f}})\, .
\end{align}
Then by \cite[Prop.8.3]{Brezis} it follows that the weak derivative $\partial_s^3 \phi$ exists and belongs to $L^2$.  By \eqref{eqa4} since $\| \Kapp\|_{\infty}, \| \partial_s \Kapp\|_{\infty}, \| \partial_s^2 \Kapp\|_{L^2} \leq C$ and $\| \phi\|_{H^2} \leq C$ we find for all  $\psi \in H^{1}_0$ integrating by parts
\begin{align}\label{eqsole}
\left| \int_{I} \partial_{s}^{2} \phi \cdot \partial_{s} \psi \diff s_{\bar{f}} \right|  \leq  \left| \int_{I} \nabla_{s}^{2} \phi \cdot \partial_{s} \psi \diff s_{\bar{f}} \right| +  C \| \psi\|_{L^2(\diff s_{\bar{f}})} \, .
\end{align}
It remains to estimate the integral on the right-hand side. 

Let $\chi:[0,1] \to [0,1]$ be smooth, $0\leq \chi\leq 1$ and such that $\chi \equiv 0$ on $[0,\frac14]$, $\chi \equiv 1$ on $[\frac34,1]$, $0 \leq \chi \leq 1$, $|\chi'|,|\chi''| \leq C $. Consider
\[\tilde{\psi}(x) := \int_0^x \psi(y) |\partial_x \bar{f}(y)| \diff y + \alpha  \chi(x) , \quad x \in [0,1]\, ,\]
with
\[ \alpha:= - \int_0^1 \psi(y) |\partial_x \bar{f}(y)| \diff y .\]
Then, using that $|\partial_x \bar{f}|\geq \delta >0$ in $I$ for some $\delta>0$,
\begin{align}\label{estpsi}
\partial_s \tilde{\psi}(x) &= \psi(x) +  \alpha\partial_s \chi(x), \, \quad | \alpha|  \leq C \| \psi\|_{L^2(\diff s_{\bar{f}})},\, \\ \nonumber
\| \tilde{\psi}\|_{L^2(\diff s_{\bar{f}})} & \leq C \| \psi\|_{L^2(\diff s_{\bar{f}})} , \quad
\|\partial_s \tilde{\psi}\|_{L^2(\diff s_{\bar{f}})} \leq C \| \psi\|_{L^2(\diff s_{\bar{f}})}\, ,
\end{align}
and since $\psi \in H^1_0$, then $\tilde{\psi} \in H^2_0$. Then we can write the integral on the right hand side of \eqref{eqsole} as
\begin{align}\label{eqthird}
 \int_{I} \nabla_{s}^{2} \phi \cdot \partial_{s} \psi \diff s_{\bar{f}} & =  \int_{I} \nabla_{s}^{2} \phi \cdot \partial_{s} (\partial_s \tilde{\psi} - \alpha \partial_s \chi) \diff s_{\bar{f}}\\  \nonumber
& = \int_{I} \nabla_{s}^{2} \phi \cdot \nabla_{s} (\nabla_s \tilde{\psi}^{\perp}  - \langle \tilde{\psi}^{\perp}, \Kapp\rangle_{\text{euc}} \partial_s \bar{f} +  \partial_s \tilde{\psi}^{\top} - \alpha \partial_s \chi) \diff s_{\bar{f}} \, .
\end{align}
Since $\tilde{\psi}^{\perp} \in H^{2,\bot}_0$ and $\phi \in D(A)$ we find
\begin{align}\label{eqpdom}
\int_{I} \nabla_{s}^{2} \phi \cdot \nabla_{s}^2 \tilde{\psi}^{\perp} \; \diff s_{\bar{f}} = a_{\bar{f}}(\phi, \tilde{\psi}^{\perp}) = \langle A \phi, \tilde{\psi}^{\perp} \rangle_{L^{2}(\diff s_{\bar{f}})} \, .
\end{align}
Writing $\psi=\psi^{\bot}+ \langle \psi, \partial_s \bar{f} \rangle_{\text{euc}} \partial_s \bar{f}$ and since $\nabla_s (\eta \partial_s \bar{f})= \eta \Kapp$ for any $\eta: I \to \R$, the other terms in \eqref{eqthird} may be written as 
\begin{align*}
& \nabla_{s} (- \langle \tilde{\psi}^{\perp}, \Kapp\rangle \partial_s \bar{f}+  \partial_s \tilde{\psi}^{\top} - \alpha \partial_s \chi) \\
 & =  2\langle \psi, \partial_s \bar{f} \rangle_{\text{euc}} \Kapp + 2\partial_s \chi \langle \alpha, \partial_s \bar{f} \rangle_{\text{euc}} \Kapp +  \langle \tilde{\psi}, \Kapp \rangle_{\text{euc}} \Kapp \\
& \quad + \langle \tilde{\psi}, \partial_s \bar{f} \rangle_{\text{euc}} \nabla_s \Kapp - \alpha^{\perp} \partial_s^2 \chi \,.
\end{align*}
From \eqref{eqthird}, since $\| \Kapp\|_{\infty}, \| \partial_s \Kapp\|_{\infty}\leq C$ with \eqref{estpsi} we obtain
\begin{align}\nonumber
\left| \int_{I} \nabla_{s}^{2} \phi \cdot \partial_{s} \psi \diff s_{\bar{f}} \right| &  \leq \| A\phi\|_{L^2(\diff s_{\bar{f}})}\| \psi\|_{L^2} + C \| \nabla_s^2 \phi\|_{L^2(\diff s_{\bar{f}})} \| \psi\|_{L^2} \\
& \leq C \| \psi\|_{L^2(\diff s_{\bar{f}})} \, .\label{eqfastfertig}
\end{align}
Combining \eqref{eqsole} and \eqref{eqfastfertig}, \eqref{ziel1} follows. Hence $D(A) \subset H^3 \cap H^{2,\bot}_0$.

To show the inclusion  $D(A) \subset H^{4, \bot}_{c}$ it remains to show that for $\phi \in D(A)$, 
\begin{align}\label{ziel2}
 \left| \int_{I} \partial_{s}^{3} \phi \cdot \partial_{s} \psi \diff s_{\bar{f}} \right|  \leq C \| \psi\|_{L^2(\diff s_{\bar{f}})} \mbox{ for all }\psi \in H^{1}_0 (\diff s_{\bar{f}}) \, .
\end{align}
At this point we use that $\bar{f} \in W^{5,2}(I;\R^d)$. By \eqref{eqa5} since $\| \Kapp\|_{\infty}, \| \partial_s \Kapp\|_{\infty}$, $\| \partial_s^2 \Kapp\|_{\infty} \leq C$ and also $\| \partial_s^3 \Kapp\|_{L^2} \leq C$, $\| \phi\|_{H^3} \leq C$ integrating by parts we find for all  $\psi \in H^{1}_0$
\begin{align}\label{eqsole2}
\left| \int_{I} \partial_{s}^{3} \phi \cdot \partial_{s} \psi \diff s_{\bar{f}} \right|  \leq  \left| \int_{I} \nabla_{s}^{3} \phi \cdot \partial_{s} \psi \diff s_{\bar{f}} \right| +  C \| \psi\|_{L^2(\diff s_{\bar{f}})} \, .
\end{align}
It remains to estimate the integral on the right hand side. 
By a density argument we can restrict to consider test-functions $\psi$ in $H^2_0$. Then
\begin{align*}
\int_{I} \nabla_{s}^{3} \phi \cdot \partial_{s} \psi \diff s_{\bar{f}}  = \int_{I} \nabla_{s}^{3} \phi \cdot \nabla_{s} \psi \diff s_{\bar{f}}  = - \int_{I} \nabla_{s}^{2} \phi \cdot \nabla_s^2 \psi \diff s_{\bar{f}} \, .
\end{align*}
Writing as before $\psi= \psi^{\top}+ \psi^{\perp}$ we compute
\begin{align*}
\nabla_s^2 \psi  & = \nabla_s^2 \psi^{\perp} + \nabla_s^2  ( \langle \psi, \partial_s \bar{f} \rangle_\text{euc} \partial_s \bar{f} )\\
& = \nabla_s^2 \psi^{\perp} + \nabla_s ( \langle \psi, \partial_s \bar{f} \rangle_\text{euc} \Kapp) 
\end{align*} 
obtaining
\begin{align*}
\int_{I} \nabla_{s}^{3} \phi \cdot \partial_{s} \psi \diff s_{\bar{f}}  & = - \int_{I} \nabla_{s}^{2} \phi \cdot \nabla_s ( \langle \psi, \partial_s \bar{f} \rangle_\text{euc} \Kapp)  \diff s_{\bar{f}}   -a_{\bar{f}}(\phi,\psi^{\perp})  \, .
\end{align*}
since $\psi^{\bot} \in H^{2,\bot}_0$. Due to the bounds on the curvature, since $\phi \in H^3$, integrating by parts and using \eqref{eq:DefAA}, it follows
\begin{align*}
\left| \int_{I} \nabla_{s}^{3} \phi \cdot \partial_{s} \psi \diff s_{\bar{f}} \right| & \leq C \| \psi\|_{L^2(\diff s_{\bar{f}})} + \| A\phi\|_{L^2(\diff s_{\bar{f}})}\| \psi\|_{L^2(\diff s_{\bar{f}})} \leq C\| \psi\|_{L^2(\diff s_{\bar{f}})} \, .
\end{align*}
Hence $D(A)=H^{4, \bot}_{c}$ and \eqref{action} holds for all $\phi \in D(A)$.

It remains to show that the norms $\| \cdot \|_{D(A)}$ and $\| \cdot \|_{H^4}$ are equivalent on $D(A)$. From its definition it is clear that $\| \phi \|_{D(A)}\leq \| \phi \|_{H^4}$ for all $\phi \in D(A)$. For the other inequality it is sufficient to show that there exists some constant $C$ such that
\[\| \partial_s^3 \phi \|_{L^2(\diff s_{\bar{f}})} + \|\partial_s^4 \phi\|_{L^2(\diff s_{\bar{f}})} \leq C( \| \phi \|_{H^2(\diff s_{\bar{f}})} +  \|\nabla_s^4 \phi\|_{L^2(\diff s_{\bar{f}})}) \, \mbox{ for all }\phi \in D(A).\]
Since $\phi$ is normal using \eqref{eqa5} and \eqref{eqa6} and the bounds on the curvature and its derivatives we find
\[\| \partial_s^3 \phi \|_{L^2(\diff s_{\bar{f}})} + \|\partial_s^4 \phi\|_{L^2(\diff s_{\bar{f}})} \leq C( \| \phi \|_{H^2(\diff s_{\bar{f}})} +  \|\nabla_s^3 \phi\|_{L^2(\diff s_{\bar{f}})}+  \|\nabla_s^4 \phi\|_{L^2(\diff s_{\bar{f}})}) \, , \]
and hence we only need that the $L^2$-norm of $\nabla_s^3 \phi$ can be controlled by $\| \phi \|_{H^2} $ and $\|\nabla_s^4 \phi\|_{L^2}$. This follows from \cite[Lemma C.4]{DP}. Indeed that result gives the existence of a constant (depending only on the length of the curve $\bar{f}$) such that for all $\varepsilon \in (0,1)$
\begin{align*} 
& \| \nabla_s^3 \phi \|_{L^2(\diff s_{\bar{f}})} \\
& \leq C [\varepsilon ( \| \nabla_s^2 \phi \|_{L^2(\diff s_{\bar{f}})} +\| \nabla_s^3 \phi \|_{L^2(\diff s_{\bar{f}})} +  \| \nabla_s^4 \phi \|_{L^2(\diff s_{\bar{f}})}) +\frac{1}{\varepsilon} \| \nabla_s^2 \phi \|_{L^2(\diff s_{\bar{f}})}] .
\end{align*}
Choosing $\varepsilon$ small enough the claim follows.
\end{proof}


\begin{proof}[Proof of Proposition \ref{prop:FredholmMain}]
By Lemmas \ref{lem2.12} and \ref{lem2.13},  the operator $\nabla_s^4+\omega :H^{4, \bot}_c(\diff s_{\bar{f}}) \to   L^{2,\bot}(\diff s_{\bar{f}})$ is an isomorphism and hence a Fredholm operator of index zero.  
Since the embedding $H^{4, \bot}_c (\diff s_{\bar{f}})\hookrightarrow L^{2, \bot}(\diff s_{\bar{f}})$ is compact 
we find that
\[
\nabla_s^4: H^{4, \bot}_c (\diff s_{\bar{f}})\to L^{2,\bot} (\diff s_{\bar{f}})
\]
is also Fredholm of index $0$ (see \cite[Example 8.16 (ii))]{Zeidler}). This yields the claim.
\end{proof}

\begin{corol}\label{cor:Fredholm}
Let $\bar{f} \in H^5$ be a regular curve. The Fr\'echet derivative 
\begin{equation}
 \mathfrak{L} := (d E)'(0):H^{4,\bot}_{c} \to L^{2,\bot}
\end{equation}
of $d E$ at zero is a Fredholm operator of index $0$.
\end{corol}
\begin{proof}
The operator $\mathfrak{L}$ is associated to the second variation $E''(0)$ (recall \eqref{ciccio}) which is given in  Proposition \ref{prop:2ndvar}. 
For $\phi \in H^{4,\bot}_{c}$ we find using partial integration and the boundary value
\begin{align}
\label{eq:Lrepresent}
\mathfrak{L} \phi & = \Big{(} \nabla_s^4 \phi 
 +  (\nabla_s^2 \phi \cdot \Kapp) \Kapp  
 +\frac{3}{2} |\Kapp|^2  \nabla_s ^2 \phi 
 + 3 (\Kapp \cdot  \nabla_s \Kapp)  \nabla_s \phi 
  \\ \nonumber & \qquad
  +2 (\nabla_s \phi \cdot \Kapp) \nabla_s\Kapp 
   + 2 (\phi \cdot \Kapp) \nabla_s^2 \Kapp  
 + 3 (\phi \cdot  \nabla_s \Kapp)  \nabla_s\Kapp 
\\ \nonumber &\qquad  
  + (\phi \cdot \nabla_s^2 \Kapp) \Kapp 
  + |\Kapp|^2 (\Kapp \cdot \phi) \Kapp  
 \Big{)} |\partial_x \bar{f}| \\ \nonumber
& = |\partial_x \bar{f}| \nabla_s^4 \phi + B(\phi) \in L^{2,\perp}\, ,
\end{align} 
for a linear operator $B: H^{4,\bot}_{c} \to L^{2,\bot}$ which is compact since the embeddings $H^4 \hookrightarrow L^2, H^1, H^2$ are all compact and the coefficients are uniformly bounded since $\bar{f} \in H^5$. By Proposition \ref{prop:FredholmMain}  and since $|\partial_x \bar{f}|$ is uniformly bounded from above and below, it follows that $H^{4,\bot}_{c} \ni \phi \mapsto  |\partial_x \bar{f}| \nabla_s^4 \phi \in L^{2,\bot}$ is a Fredholm operator of index $0$. Since $B$ is compact, the claim follows using that the sum of a Fredholm operator of index zero and a compact operator is again 
a Fredholm operator of index $0$.
\end{proof}


\section{Proof of the \L ojasiewicz-Simon inequality for all directions}
\label{sec4}

In the previous section we have shown the \L ojasiewicz-Simon inequality for the functional $E$ (see \eqref{Esimple}, Theorem \ref{LosNormal} and Remark \ref{remgeometric}), i.e.\ we have considered only variations in the normal direction. This is needed to get 
the desired Fredholm property of the second variation.  In this section we want to show the existence of constants $C_{1} \geq 0, \theta \in (0,\frac{1}{2}]$ and $\sigma > 0$ such that the \L ojasiewicz-Simon inequality
\begin{equation} 
 \betrag{\mathbf{E}(0)-\mathbf{E}(\psi)}^{1-\theta} \leq C_{1} \norm[L^2(\diff s_{\bar{f} + \psi})]{\nabla \mathbf{E}(\psi)}. \tag{\ref{eq:LS}}
\end{equation}
is actually satisfied on a $\sigma$-ball around zero of the \emph{whole} space of variations 
$\psi \in H^{4}_\text{c}$. 
This can be achieved starting from Theorem \ref{LosNormal} and Remark \ref{remgeometric} and by noticing that variation vector fields that are tangent to a fixed immersion $\bar{f}$ correspond to reparametrizations of $\bar{f}$. 

\begin{lemma} \label{lemma:reparametrization}
 Let $\bar{f} \in H^5(I;\R^d)$ be a regular curve. There exists a $\sigma= \sigma (\bar{f})>0$ such that for any $\psi \in H^4_\text{c}$ with $\norm[H^4]{\psi} \leq \sigma$, there exists a $H^4$-diffeomorphism $\Phi:I\to I$ such that 
 \begin{equation}\label{repara}
  (\bar{f}+\psi)\circ \Phi = \bar{f}+\phi
 \end{equation}
for some $\phi \in H^{4,\bot}_{\text c}$. 

Moreover, for $\tilde \sigma>0$  given, there exists $\sigma= \sigma (\tilde \sigma,\bar{f})>0$ such that for any $\psi \in H^4_\text{c}$ with $\norm[H^4]{\psi} \leq \sigma$, \eqref{repara} is valid and the normal vector field $\phi$ satisfies the inequality $\norm[H^4]{\phi} \leq \tilde \sigma$.
\end{lemma}

\begin{remark}\label{remarkDiff}
As can be easily seen from the proof (see \eqref{auf1} below), one can achieve higher regularity of the diffeomorphism in the previous lemma by assuming more regularity in the data.  More precisely, if for some $m \geq 4$, $m\in \N$, $\bar{f} \in H^{m+1}(I;\R^d)$ then there exists a $\sigma= \sigma (\bar{f})$ such that for any $\psi \in H^4_\text{c}\cap H^m$ with $\norm[H^4]{\psi} \leq \sigma$, there exists a $H^m$-diffeomorphism $\Phi$ and a $\phi \in H^{4,\bot}_{\text c}\cap H^m$ so that \eqref{repara} is valid. 
Similarly (adapting the arguments given in \ref{app:proofs}) also the second part of the claim of Lemma \ref{lemma:reparametrization} remains true, that is for given $\tilde \sigma>0$ there exists $\sigma= \sigma (\tilde \sigma,\bar{f})>0$ such that for any $\psi \in H^4_\text{c} \cap H^m$ with $\norm[H^m]{\psi} \leq \sigma$, \eqref{repara} is valid and the normal vector field $\phi$ satisfies the inequality $\norm[H^m]{\phi} \leq \tilde \sigma$.
\end{remark}

To show the existence of $\Phi$, we use the implicit function theorem in the following form.\begin{mythm}[{\cite[Theorem 4.B]{Zeidler}}]\label{thmZeid}
 Let $X, Y, Z$ be real Banach spaces, $(x_0, y_0) \in X \times Y$, $\Lambda \times \Omega$ be an open neighbourhood of $(x_0, y_0)$ in $X \times Y$ and $F:\Lambda \times \Omega \to Z$ such that
\begin{enumerate}
\item $F(x_0, y_0) = 0$;
 \item $\partial_y F$ exists as partial Fr\'echet-derivative on $\Lambda \times \Omega$ and $\partial_y F(x_0,y_0): Y \to Z$ is bijective;
 \item $F$ and $\partial_y F$ are continuous at $(x_0, y_0$).
\end{enumerate}
Then there exist positive numbers $r_0$ and $r$ such that $B_{r_0}(x_0)\times B_r(y_0) \subset \Lambda \times \Omega$ and for every $x \in X$ satisfying $\norm[X]{x-x_0} \leq r_0$, there is exactly one $y=y(x) \in Y$ for which $\norm[Y]{y-y_0} \leq r$ and $F(x,y) = 0$. Moreover, if $F$ is continuous in a neighborhood of $(x_0,y_0)$, then $y(\cdot)$ is continuous in a neighborhood of $x_0$.
\end{mythm}

\begin{proof}[Proof of Lemma \ref{lemma:reparametrization}]
Let $X := H^4_{\text{c}}$, $Y := Z := H^3\cap H^1_0(I;\R)$ and $(x_0, y_0) := (0,0) \in X \times Y$. Let $\Lambda := B_\rho(0)\subset X$ be small enough such that $\bar{f}+\psi$ is immersed for all $\psi \in \Lambda$ and $\Omega := B_R(0) \subset Y$ be small enough such that $\norm[\infty]{\varphi'} < \frac{1}{2}$ for all $\varphi \in \Omega$. Note that this is possible since $Y \hookrightarrow C^1(I;\R)$. Moreover, notice that this choice of $\Omega$ implies that $\id_I+\varphi$ is a $C^1$-diffeomorphism of $I$ for all $\varphi \in \Omega$.

Consider the functional
\begin{equation}
F:\Lambda\times\Omega \to Z,\;\; (\psi,\varphi) \mapsto \metric[\text{euc}]{(\bar{f}+\psi)\circ (\id_I + \varphi) - \bar{f}, \partial_x \bar{f}},
\end{equation}
which is well defined since the composition of a functions $H^m(I,\R)$ with a function in $H^{m}$, that is also a $C^{1}$-diffeomorphism, is an element of $H^{m}(I,\R)$ for $m\geq 2$. Indeed, by \cite[Prop.9.5]{Brezis}   
we see that since $\bar{f}+\psi \in H^1$, we have $(\bar{f}+\psi)\circ (\id_I +\varphi) \in H^1$ and its weak derivative is (as expected) given by
\[(1+\varphi') \, (\bar{f}+\psi)'\circ (\id_I +\varphi) \, .\]
Since $(\bar{f}+\psi)' \in H^1$ we can repeat the same argument and we obtain $(\bar{f}+\psi)\circ (\id_I +\varphi) \in H^2$. The case $m>2$ can be treated repeating the same arguments. From the definition, one sees also that the function $F$ is continuous.

We now show that with these choices the assumptions of Theorem \ref{thmZeid} are satisfied. By definition $F(0,0) = 0$ and writing for $h$ sufficiently small
\begin{align*}
& \quad F(\psi, \varphi + h)- F(\psi , \varphi) -  h   \metric[\text{euc}]{(\partial_x (\bar{f}+\psi))\circ (\id_I + \varphi)  , \partial_x \bar{f}}\\
& = \sum_{i=1}^n [\int_{\id_I+\varphi}^{\id_I+\varphi+h} \partial_x^2 ((\bar{f}+\psi)^i)(t) (\id_I +\varphi +h -t) \; \diff t ] \partial_x \bar{f}^i
\end{align*}
we see that the Fr\'echet derivative of $F$ with respect to the second component exists and is given by
\begin{align*}
\partial_y F(\psi,\varphi): &  H^{3}(I,\R)\cap H^{1}_0(I,\R)  \to H^3\cap H^1_0(I;\R) \\
&h \mapsto h   \metric[\text{euc}]{(\partial_x (\bar{f}+\psi))\circ (\id_I + \varphi)  , \partial_x \bar{f}}\, ,
\end{align*}
a linear and continuous operator. In particular, $\partial_y F(0,0)$ is the scalar multiplication with $|\partial_x \bar{f}|^2$ acting from $H^3\cap H^1_0$ to $H^3\cap H^1_0$. This is an invertible and continuous operator since $\bar{f}$ is an immersion.

Since $F$ and $F_{y}$ are continuous in a neighborhood of $(0,0)$, it follows from Theorem \ref{thmZeid} that there exist some $0 < r_0 \leq \rho$ and  $0 < r \leq R$ such that for any $\psi \in H^4_{c}$ with $\|\psi\|_{H^4} \leq r_0$ there exists exactly one $\tilde\varphi: I \to \R \in H^3\cap H^1_0$, $\norm[\infty]{\tilde\varphi'} < \frac{1}{2}$ and $\| \tilde{\varphi}\|_{H^3} \leq r$, such that
\begin{equation}\label{eq:FNullstelle}
 \metric[\text{euc}]{(\bar{f}+\psi)\circ (\id_I + \tilde\varphi) - \bar{f}, \partial_x \bar{f}} = 0 \text{ on } I.
\end{equation}
Moreover, using the continuity of $\partial_y F$ we may choose $r,r_0$ small enough such that
\begin{equation}
 \metric[\text{euc}]{\partial_x (\bar{f}+\psi)\circ (\id_I + \tilde\varphi) , \partial_x \bar{f}} \geq \frac12 |\partial_x \bar{f}|^2 \ne  0 \text{ on } I.
\end{equation}
We have already used that $\id_I + \tilde{\varphi}$ is a diffeomorphism on $I$. Since $F$ is continuous, $\| \tilde{\varphi}\|_{H^3}$ depends continuously on $\| \psi\|_{H^4}$. We show now that $\tilde{\varphi}$ is actually in $H^4$ and that also the $H^4$-norm of $\tilde{\varphi}$ depends continuously on $\| \psi\|_{H^4}$.

Differentiating \eqref{eq:FNullstelle} one sees that
\begin{align}\label{auf1}
1+\tilde{\varphi}'   =\frac{ |\partial_x \bar{f}|^2 -  \metric[\text{euc}]{(\bar{f}+\psi)\circ (\id_I + \tilde\varphi) - \bar{f}, \partial_x^2 \bar{f}}}{\metric[\text{euc}]{\partial_x(\bar{f}+\psi)\circ (\id_I + \tilde\varphi) , \partial_x \bar{f}}} \, .
\end{align}
Since the right-hand side is in $H^3(I,\R)$, it follows that $\tilde{\varphi} \in H^4(I,\R)$.
The first part of the claim follows by letting $\Phi := \id+\tilde\varphi$ and $\phi := (\bar{f}+\psi)\circ \Phi - \bar{f}$. Then $\Phi$ is a $H^4$-diffeomorphism and, by construction (see \eqref{eq:FNullstelle}) $\phi$ is normal to $\partial_x \bar{f}$. 

The proof of the second part of the claim is quite technical and it is given in the appendix~\ref{app:proofs}.
\end{proof}

By virtue of this lemma we may now derive the  \L ojasiewicz-Simon inequality for all directions from the one already proven for normal directions. 

\begin{proof}[Proof of Theorem \ref{mainLos}]
As already observed, since $\bar{f} \in H^4(I; \R^d)$ is a regular curve and a critical point of the elastic energy, it is smooth.    
By Theorem \ref{LosNormal} there exists a $\tilde{\sigma}>0$ and constants $C_2, \theta$ such that for all $\phi \in B_{\tilde{\sigma}}(0)\subset H^{4,\bot}_{c}$ inequality  \eqref{eq:LSNormal} holds. By choosing $\sigma$ as in  Lemma \ref{lemma:reparametrization} we find that for any 
$\psi \in B_{ \sigma}(0) \subset H^4_\text{c}$  there exists a $H^4$-diffeomorphism $\Phi$ as well as a normal vector field $\phi$ with $\| \phi\|_{H^4} \leq \tilde{\sigma}$ such that, by the geometric invariance of the functional $\mathcal{E}$, we can write
\begin{align*}
 \mathbf{E}(\psi) &= \mathcal{E}(\bar{f}+\psi) = \mathcal{E}((\bar{f}+\psi)\circ \Phi) = \mathcal{E}(\bar{f}+\phi)= E(\phi). 
\end{align*}
Hence, with the constants $C_2$, $C_{3}$ and $\theta$ from Theorem \ref{LosNormal} and Remark~\ref{remgeometric} we obtain
\begin{align*}
|  \mathbf{E}(\psi)  - \mathbf{E}(0)|^{1-\theta}   =  | E(\phi)  -E(0)|^{1-\theta} 
&\leq C_2 \norm[L^2(dx)]{d E (\phi)}  \leq C_{3} \| \nabla E (\phi) \|_{L^{2}(\diff_{\bar{f} +\phi })}.
 \end{align*}
Recalling the explicit formulas for the gradients given below Definition~\ref{def2}  we see that $ \| \nabla E (\phi) \|_{L^{2}(\diff s_{\bar{f} +\phi })} \leq \| \nabla \mathbf{E} (\phi) \|_{L^{2}(\diff s_{\bar{f} +\phi })}$.
Using now the geometric invariance of the gradient of the energy we find that $  \| \nabla \mathbf{E} (\phi) \|_{L^{2}(\diff s_{\bar{f} +\phi })} = \norm[L^2(\diff s_{\bar{f}+\psi})]{\nabla \mathbf{E} (\psi)}$ and therefore
\begin{align*}
|  \mathbf{E}(\psi)  - \mathbf{E}(0)|^{1-\theta}  & \leq C_3 \norm[L^2(\diff s_{\bar{f}+\psi})]{\nabla \mathbf{E} (\psi)} \,.
\end{align*} 
The claim follows with these choices of $\sigma$ and $\theta$ and choosing $C_1= C_3$.
\end{proof}


\section{Main result: convergence to elastica}\label{sec5}
In this section we turn to the geometric problem under consideration, that is we consider the evolution of smooth regular open curves with fixed endpoints, fixed unit tangents at the boundary, and moving according to the $L^2$-gradient for the elastic energy $\E_{\lambda}$. A precise formulation has already been given in \eqref{eq:ElasticFlow}.
  Note that now the arc-length element $\diff s_f =|f_x (x,t)| \diff x$ is time dependent. 
Moreover note that the claim of Theorem~2.4 holds also for the functional $\E_\lambda$ since only  lower order terms have been added (recall Corollary~\ref{cor:2varl}).

By Theorem~\ref{ThmLin} one is able to find sequences of curves $f(t_i)$, $t_i \to \infty$, converging smoothly (after an appropriate reparametrization) to  a smooth regular curve $f_\infty$.
As we mentioned in the introduction, in principle for different (sub-)sequences we could find different limits. Here we want to show that this is not the case: that is, for a chosen initial data $f_0$ the flow evolves (after a suitable reparametrization) to one critical point. This claim can be achieved by application of the \L ojasiewicz-Simon inequality as described in the following.

\begin{proof}[Proof of Theorem~\ref{main2}]
In the proof of long-time existence of the flow it is shown that  the length of the curve remains uniformly bounded from above and below along the evolution (with constants independent of time: see \cite[equations (41), (43), (47)]{Lin}).
However, the arc-length element $|f_x|$ might degenerate in the limit and this is the reason why a reparametrization of  the curves $f(t_i)$ is necessary for the sub-convergence result.
 (In fact uniform bounds in time for $f$ and its derivatives are obtained when the curve is considered reparametrized by arc-length.)
With this in mind let us introduce the smooth map $\Psi: [0,\infty) \times [0,1] \to [0,1]$,
\begin{align*}
\Psi(t,x):=\frac{1}{\mathcal{L}(f(t))} \int_0^x |f_x (t, \xi)| d\xi.
\end{align*}
For any fixed $t$ the map $\Psi$ gives a smooth diffeomorphism by which we can reparametrize the curve $f(t)$.\\
As already pointed out, by Theorem \ref{ThmLin} we know that the flow exists for all times and it is smooth.
By construction of the $L^{2}$-gradient flow we have that
\begin{align}
\frac{\diff}{\diff t}\E_\lambda(f(t)) &= \langle \nabla_{L^{2}} \E_\lambda(f(t)), \partial_t f \rangle_{L^{2}(\diff s_{f})} 
= - \norm[L^2 (\diff s_{f})]{\partial_t f}^2 \leq 0. \label{eq:monotone}
\end{align}
Thus $\E_\lambda(f(t)) $ decreases in time. Moreover by the subconvergence result there exists a sequence $t_i \to \infty$ such that 
\begin{equation}\label{subconvergence}
f(t_i,\Psi^{-1}(t_i, \cdot) ) \mbox{ converges to }f_{\infty} \mbox{ in } C^k \mbox{ for all }k \in \N,
\end{equation}
with $f_{\infty}$ a critical point  of the functional $\mathcal{E}_{\lambda}$. (This critical point $f_{\infty}$ will be kept fixed throughout the following arguments.) It then follows that
\begin{equation}\label{behenergy} 
\mathcal{E}_{\lambda} (f(t)) = \mathcal{E}_{\lambda} (f(t, \Psi^{-1}(t, \cdot))) \geq   \mathcal{E}_{\lambda} (f_{\infty})  \mbox{ for all }t\, . 
\end{equation}
In particular we observe that possibly different critical points, that are limit to different sequences of times, share the same energy level $\mathcal{E}_{\lambda} (f_{\infty})$. We distinguish now between two cases, namely whether the final energy level is attained in finite time or not.

If there exists $\bar{t} \in (0,\infty)$ such that  $\mathcal{E}_{\lambda} (f(\bar{t})) =  \mathcal{E}_{\lambda} (f_{\infty})$, then necessarily  $\mathcal{E}_{\lambda} (f(t)) =  \mathcal{E}_{\lambda} (f_{\infty})$ for all $t \in [\bar{t},\infty)$. Then $f(t,\cdot)=f(\bar{t}, \cdot)$ for all $t \in [\bar{t},\infty)$ since $\partial_t f \equiv 0$. Due to the subconvergence result we additionally find
\[ f(t, \Psi^{-1}(t,\cdot))= f_{\infty} \mbox{ for all }t \geq \bar{t}, \]
and the claim follows. 

Hence in the following 
we may assume that  $\mathcal{E}_{\lambda} (f(t)) > \mathcal{E}_{\lambda} (f_{\infty})$ for all $t$. By Theorem~\ref{mainLos} (now adapted with the obvious changes to the functional $\E_\lambda$) and since $f_\infty$ is a regular critical point for $\E_\lambda$ we conclude that
there exists $\theta \in (0,\frac{1}{2}]$ and $C_1,\sigma > 0$ such that the \L ojasiewicz-Simon gradient inequality 
\begin{equation}\label{6.7}
 \betrag{\mathbf{E}_\lambda(0)-\mathbf{E}_\lambda(\psi)}^{1-\theta} \leq C_1 \norm[L^2 (\diff s_{(f_\infty+ \psi)} )]{\nabla \mathbf{E}_{\lambda}(\psi)}
\end{equation}
holds for all $\psi \in H^{4}_\text{c}$ such that $\norm[H^4]{\psi}\leq \sigma$. Here $\mathbf{E}_\lambda(\psi)= \E_{\lambda} (f_{\infty} +\psi)= \E (f_{\infty} +\psi) +\lambda \mathcal{L}(f_{\infty} +\psi)$ and \eqref{6.7} can be written as \begin{equation}\label{6.8}
 \betrag{\E_{\lambda} (f_{\infty})-\E_{\lambda} (f_{\infty} +\psi)}^{1-\theta} \leq C_1 \norm[L^2 (\diff s_{(f_\infty+ \psi)} )]{\nabla_{L^{2}} \E_{\lambda} (f_{\infty} +\psi)}.
\end{equation}

For reasons that will become clearer shortly and in order to highlight the dependence of some important parameters, let us now fix $\delta=\delta(f_{\infty}) >0$, such that for any map $\tilde{f}$ with
\[ \| \tilde{f} - f_{\infty} \|_{C^{1}[0,1]}  < \delta\]
we have that $|\partial_{x} \tilde{f}| \geq M > 0$  with $M := \frac{1}{2} \min_{I}|\partial_{x} f_{\infty}|$ (recall that $f_{\infty}$ is a regular curve), as well as 
\begin{align}\label{nup}
| \nu^{\perp} | & \geq \frac{1}{2} |\nu| \, 
\end{align}
for any $\nu$ vector field that is normal to $f_{\infty}$.
Here $\nu^{\perp}$ denotes the normal component of $\nu$ with respect to $\tilde{f}$.
 Note that the above inequality makes sense, since  for $\nu$ orthogonal to $\tau_{\infty}:=\frac{\partial_x f_{\infty}}{| \partial_x f_{\infty} |}$, we find
\begin{align*}
\nu^{\perp} & = \nu - \langle \nu, \partial_s \tilde{f} \rangle \partial_s \tilde{f}  = \nu - \langle \nu, \partial_s \tilde{f} - \tau_{\infty} \rangle \partial_s \tilde{f},
\end{align*} 
and $| \nu^{\perp} |  \geq (1-a) |\nu| $
for $a$ small when $\| \tilde{f}- f_{\infty}\|_{C^1}$ is small enough.
Next, without loss of generality we may assume that   $\sigma < \delta$, with $\sigma$ as in the 
\L ojasiewicz-Simon gradient inequality \eqref{6.8} above.

Let $\varepsilon>0$ to be chosen. From \eqref{subconvergence} it follows that there exists $k \in \mathbb{N}$ such that 
\begin{equation}\label{startnear}
\| f(t_k, \Psi^{-1}(t_k,\cdot))- f_{\infty}\|_{H^{5}} < \varepsilon \, .
\end{equation}
With this fixed choice of $k=k(\varepsilon)$, we define
\[ \hat{f}_0 (x):= f(t_k, \Psi^{-1}(t_k,x)) \mbox{ for }x \in [0,1] \, .\]
Then $ \hat{f}_0$ is a regular smooth curve satisfying the clamped boundary conditions of~\eqref{eq:ElasticFlow}. For the next steps, we follow quite closely the method of proof presented in \cite{blowups}. 
The idea is to start with $\hat{f}_{0}$ as a initial data, and show the existence of  an elastic flow  that can be written as a graph over $f_{\infty}$. This is possible because we start really close to the critical point (closeness in norm and hence, so to say, in parametrization). Subsequently we will show that this flow exists for all time and differs from the original one by a suitable reparametrization. The reason for this somehow cumbersome ansatz lies in the fact that although we already have long-time existence of the flow for $f$, we have very little control on  its parametrization. In particular, once we fix a limit point $f_{\infty}$, we automatically pick a parametrization and there is no reason to believe that $f$ should converge to it in suitable norms.   

By Lemma \ref{lemma:reparametrization} and Remark \ref{remarkDiff} and choosing $\varepsilon < \sigma(f_{\infty})$ (with $\sigma(f_\infty)$ as defined in Remark \ref{remarkDiff} with $m=5$) there exists a diffeomorphism $\Phi_0 \in H^5(I)$ such that 
\[\hat{f}_0 \circ \Phi_0 = f_{\infty} + N_0\]
with $N_0 \in H^5(I;\R^d)$ and normal along $f_\infty$. Moreover, the same lemma and remark give the existence of a constant $K(\varepsilon)$ such that $K(\varepsilon) \searrow 0$ for $\varepsilon \searrow 0$ and such that
\begin{equation}\label{N0est}
\| N_0\|_{H^5} = \| \hat{f}_0 \circ \Phi_0 - f_{\infty} \|_{H^5} < K(\varepsilon)< \frac{1}{2} \sigma <\frac{1}{2}\delta\, 
\end{equation} 
for $\varepsilon$ small enough.
Since $f_{\infty}$ is smooth and also $\hat{f}_0$ is smooth, proceeding as in the proof of Lemma~\ref{lemma:reparametrization} one sees that $\Phi_0$ and $N_0$ are also smooth.

With \eqref{eq:ElasticFlow} in mind, we now look for a solution $\tilde{f}$ of
\begin{equation}\label{solnormal}
\left\{ \begin{array}{ll}
\partial_t^{\perp} \tilde{f} = - \nabla_{L^{2}} \mathcal{E}_{\lambda} (\tilde{f}) & \mbox{ in } (0,T)\times I,\\
\tilde{f}(0,\cdot) = \hat{f}_0\circ\Phi_0 = f_{\infty}+N_0  & \mbox{ on } I\\
\tilde{f}(t,0)=f_{-}, \tilde{f}(t,1)=f_{+} & \mbox{ for } t \in (0,T)\\
\partial_s \tilde{f}(t,0)=T_{-}, \partial_s \tilde{f}(t,1)=T^{+}& \mbox{ for } t \in (0,T)
\end{array}
\right. 
\end{equation}
of  type $\tilde{f}= f_{\infty}+ N$ with $N$ normal to $f_{\infty}$. Here $\partial_t^{\perp} \tilde{f} := \partial_t \tilde{f} - \langle \partial_t \tilde{f}, \partial_s \tilde{f}\rangle \partial_s \tilde{f} $
 with $s=s_{\tilde{f}}$. Moreover $T_{-}=\frac{(f_{\infty})_{x}}{|(f_{\infty})_{x}|}(0)$, $T_{+}=\frac{(f_{\infty})_{x}}{|(f_{\infty})_{x}|}(1)$, $f_{\infty}(0)= f_{-}$, $f_{\infty}(1)=f_{+}$. 

Let $\{ N_i\}_{i=1}^{d-1}$ be a smooth (time-independent) frame in $\R^d$ orthogonal to $f_{\infty}$ and orthonormal. Then we may write
\[ \tilde{f} =f_{\infty}+ N= f_{\infty} + \sum_{i=1}^{d-1} \varphi^i N_i \mbox{ with } \varphi^i: [0,T) \times I \to \R \, .\]
We find
\[
 \partial_t \tilde{f}   =  \sum_{i=1}^{d-1} (\partial_t \varphi^i) N_i , \quad \partial_t^{\perp} \tilde{f}  =  \sum_{i=1}^{d-1} (\partial_t \varphi^i) N_i^{\perp}
\]
and there exists some vector-valued function $\vec{P}$ depending smoothly on its arguments as well as $f_\infty$ and the frame $\{ N_i\}_{i=1}^{d-1}$ such that
\[
\nabla_{L^{2}} \mathcal{E}_{\lambda} (\tilde{f})  = \frac{1}{|\partial_x (f_{\infty}+N)|^4} \sum_{i=1}^{d-1} (\partial_x^4 \varphi^i) N_i^{\perp} + \vec{P}(\cdot, \varphi, \partial_x \varphi, \ldots, \partial_x^3 \varphi) \, ,
\]
where we abbreviate $\varphi = (\varphi^1,\ldots, \varphi^{d-1})$.

We claim that for $t$ sufficiently small, $\{ N_i^{\perp}\}_{i=1}^{d-1}$ is a basis of $\{ \partial_s \tilde{f}\}^{\perp}$. Indeed we have seen that \eqref{nup} holds when $\| \tilde{f}- f_{\infty}\|_{C^1}$ is small enough. 
Hence in a short interval of time $\{ N_i^{\perp}\}_{i=1}^{d-1}$ is a basis of the normal bundle of $\tilde {f}$
and \eqref{solnormal} becomes equivalent to the following PDE-system for $\varphi =(\varphi^1,\ldots, \varphi^{d-1})$:
\begin{equation}\label{system}
\left\{ \begin{array}{ll}
\partial_t \varphi = - \frac{1}{|\partial_x (f_{\infty}+N)|^4} \partial_x^4 \varphi- P(\cdot, \varphi, \partial_x \varphi, \ldots, \partial_x^3 \varphi) & \mbox{ in } (0,T)\times I,\\
\varphi(0,\cdot) = \varphi_0  & \mbox{ on } \{0\}\times I\\
\varphi(t,0)=0=\varphi(t,1) & \mbox{ for } t \in (0,T)\\
\partial_x \varphi(t,0)=0=\partial_x \varphi(t,1)& \mbox{ for } t \in (0,T)
\end{array}
\right. 
\end{equation} 
with some  function $P$ with the same dependencies as $\vec{P}$ above and with  
 $N_0= \sum_{i=1}^{d-1} \varphi_0^i N_i $. Since $|\partial_x (f_{\infty}+N)|$ is bounded away from zero  for  $\|N\|_{C^{1}}$  small enough (which is guaranteed at the initial time by \eqref{N0est}) and since the initial data belongs to $C^{4, \frac{1}{2}}([0,1])$ by embedding results, one can show that there exists some maximal time interval $(0, T')$, $0< T' \leq \infty$, such that the above system admits a unique solution in the parabolic H\"older space $C^{1+\alpha/4, 4+ \alpha}([0,T')\times[0,1])$ for $\alpha \in(0,\frac12]$ (notation as in  \cite[\S~5.1]{Lunardi}). This short-time existence result follows from classical results for strictly  parabolic quasi-linear systems (\cite{Lunardi}, \cite{lady}, \cite{eidelman}). The main strategy involves obtaining optimal regularity results for the linear parabolic problem and then applying a fixed-point argument in the appropriate H\"older spaces. A detailed proof in the context of Sobolev-spaces will soon appear in the PhD-Thesis of Spener.

Let $T\in (0,T']$ be maximal such that
\begin{equation}\label{defT}
\| \tilde{f}(t)-f_{\infty}\|_{H^4}   < \sigma < \delta \,  \mbox{ for all }t \in [0,T).
\end{equation}
In particular observe that the above inequality implies that $\{ N_i^{\perp}\}_{i=1}^{d-1}$ is a basis of $\{ \partial_s \tilde{f}\}^{\perp}$ since $\sigma < \delta$ (recall \eqref{nup}), as well as $|\tilde{f}(t)_{x}| \geq M>0$ uniformly in time.

We want to show that $T=T'=\infty$. Assume that this is not the case, hence either $0<T< T'=\infty$ or $0<T \leq T' < \infty$. Let $T'':= \min\{T', T+1\}$. 
We claim now that 
\begin{equation}\label{514bb} \sup_{[0,T)} \| \varphi(t) \|_{C^{4,\gamma}([0,1])} \leq C(\sigma, T''), 
\end{equation} and hence
\begin{align}\label{514b}
 \sup_{[0,T)} \| \tilde{f}(t) - f_{\infty} \|_{C^{4,\gamma}([0,1])} \leq C(\sigma, T'') 
 \end{align}
for some $\gamma \in (0,\frac12)$. The above upper bounds can be obtained by applying parabolic Schauder estimates. For completeness, we present here a possible method to derive such estimates. First of all note that  
by \eqref{defT} and embedding theory the map $\tilde{f}$ and hence $\varphi$ belong to  $C^{3,\frac{1}{2}}([0,1])$ uniformly in time on $[0,T]$.
In particular $\sup_{[0,T]}\|P(\cdot, \varphi, \partial_x \varphi, \ldots, \partial_x^3 \varphi) \|_{C^{0,\frac{1}{2}}([0,1])} \leq C\sigma$ and $\sup_{[0,T]} \| \frac{1}{|\partial_x (f_{\infty}+N)|^4} \|_{C^{0,\frac{1}{2}} ([0,1])} \leq C \sigma$.  By \cite[Theorem~2.1 with $g(t,x) =\varphi_0(x)$ for all $t$]{DongZ}  and \eqref{N0est}
it follows that  $\varphi$ satisfies (in the notation of \cite{DongZ})
\begin{align*}
[\varphi]_{\frac{2+\frac{1}{2}}{4}, \frac{1}{2}, [0,T]\times [0,1]}+ [\varphi_{x}]_{\frac{1+\frac{1}{2}}{4}, \frac{1}{2}, [0,T]\times [0,1]} + [\varphi_{xx}]_{\frac{\frac{1}{2}}{4}, \frac{1}{2}, [0,T]\times [0,1]}  
&  \leq C(\sigma, T'').
\end{align*}
That is, the derivatives with respect to $x$ up to order two of $\varphi$ are H\"older continuous in time and space. In order to apply standard regularity theory we need the same regularity result also for the third derivative with respect to $x$ of $\varphi$. We obtain this as follows. Since $\partial_x \varphi$ is continuous in time and space, by \cite[Thm. 2.1]{DHP} we obtain for all $p\in (1,\infty)$
\[ \varphi \in W^{1,p}(0,T; L^p(0,1)) \cap L^{p}(0,T; W^{4,p}(0,1))\, .\]
By real interpolation (see \cite[Prop.~1.1.3]{LunInter}, \cite[Thm.~3.1]{Amann00} and \cite[2.3.1(7)]{Triebel78}) one has for $\theta \in (0,1)$ such that $4\theta \not \in\N$
\begin{align*} 
&  W^{1,p}(0,T; L^p(0,1)) \cap L^{p}(0,T; W^{4,p}(0,1)) \\
& \hookrightarrow (W^{1,p}(0,T; L^p(0,1)) , L^{p}(0,T; W^{4,p}(0,1)) )_{\theta,p}\\ & \qquad = W^{1-\theta,p}(0,T; W^{4\theta,p}(0,1)) \, .
\end{align*} 
Letting $p > 5$ we find some $\theta \in (0,1)$ such that the Sobolev Embedding Theorem \cite[(3.2),(3.3)]{Amann00} yields
\[ W^{1-\theta,p}(0,T; W^{4\theta,p}(0,1))  \hookrightarrow C^{\gamma_1}([0,T]; C^{3+\gamma_2}[0,1])\]
for some $\gamma_1,\gamma_2 \in (0,\frac12)$.
We have finally obtained that the third derivatives of $\varphi$ are H\"older in time and space on $[0,T] \times [0,1]$. Hence we can now see the PDE-system \eqref{system} as a linear system with H\"older coefficients. Hence by classical results as \cite[Theorem VI.21]{eidelman} or  \cite[Thm. 4.9]{Sol} we get \eqref{514bb} and hence \eqref{514b}.
 
We may now finally employ the \L ojasiewicz-Simon gradient inequality to finish the proof. To do so we let
\[G(t) := \left(\E_\lambda(\tilde{f}(t))-\E_\lambda(f_\infty)\right)^\theta \, \mbox{ for }t \in (0,T)\]
with $\theta$ the parameter from the  \L ojasiewicz-Simon gradient inequality (see \eqref{6.7}). Since $\tilde{f}$ is in fact a reparametrization of $f$ as we will see below after \eqref{limitevero} we may infer from \eqref{eq:monotone} and \eqref{behenergy} that $G$ is monotonically decreasing to zero.
Moreover for $t \in (0, T)$ using \eqref{6.8},  \eqref{solnormal}, and \eqref{defT}
we get
\allowdisplaybreaks{\begin{align}
 -\frac{\diff}{\diff t} G(t) &  = - \theta G(t)^{\frac{\theta-1}{\theta}}\metric[L^2 (\diff s_{\tilde{f}})]{\nabla_{L^2 } \E_\lambda(\tilde{f}(t)),\partial_t \tilde{f}} \nonumber \\
&  = - \theta G(t)^{\frac{\theta-1}{\theta}}\metric[L^2 (\diff s_{\tilde{f}})]{\nabla_{L^2 } \E_\lambda(\tilde{f}(t)),(\partial_t \tilde{f})^\perp }
   \nonumber \\
& = \theta G(t)^{\frac{\theta-1}{\theta}}\norm[L^2 (\diff s_{\tilde{f}})]{\nabla_{L^2} \E_\lambda(\tilde{f}(t))}\norm[L^2 (\diff s_{\tilde{f}})]{(\partial_t \tilde{f})^\perp} \nonumber \\
  & \geq \theta G(t)^{\frac{\theta-1}{\theta}}\frac{1}{C_1} \left(\E_\lambda(\tilde{f}(t)) - \E_\lambda(f_\infty)\right)^{1-\theta} 
  \norm[L^2 (\diff s_{\tilde{f}})]{(\partial_t \tilde{f})^\perp} 
\nonumber   \\
  &= \frac{\theta}{C_1} \norm[L^2 (\diff s_{\tilde{f}})]{(\partial_t \tilde{f})^\perp} . \label{eq:ineqG}
   \end{align}}  
Since $\langle N, \tau_{\infty}\rangle =0$, then also $\langle \partial_t N, \tau_{\infty}\rangle =0$ and hence
\begin{align*}
(\partial_t \tilde{f})^\perp & =\partial_t  N - \langle \partial_t N, \partial_s \tilde{f} \rangle  \partial_s \tilde{f} \\
& = \partial_t  N - \langle \partial_t N, \partial_s \tilde{f} - \tau_{\infty} \rangle  \partial_s \tilde{f} \, .
\end{align*}
By \eqref{nup} and \eqref{defT} we get
\[\norm[L^2 (\diff s_{\tilde{f}})]{(\partial_t \tilde{f})^\perp} \geq \frac12 \norm[L^2 (\diff s_{\tilde{f}})]{\partial_t N} = \frac12\norm[L^2 (\diff s_{\tilde{f}})]{\partial_t \tilde{f}}\, . \]
Hence we find by \eqref{eq:ineqG} and the boundedness of $|\tilde{f}_{x}(t)|$ from below 
\begin{equation}\label{derGest}
-\frac{\diff}{\diff t} G(t)  \geq \frac{\theta \sqrt{M}}{2 C_1 }
\| \partial_t \tilde{f}\|_{L^{2}} \mbox{ for }t \in (0,T) \, . 
\end{equation}
This inequality together with an interpolation argument allows us to prove that $T=\infty$. 
We first observe that for $t \in (0,T)$ by \eqref{derGest} and the definition of $G$ \begin{align}\nonumber
 \norm[L^2]{\tilde{f}(t)- f_\infty} 
&\leq \norm[L^2]{\tilde{f}(0)- f_\infty} + \int_0^t   \norm[L^2]{\partial_t\tilde{f}} \diff \tau \\ \nonumber
& \leq  \norm[L^2]{\tilde{f}(0)- f_\infty} - C (G(t)-G(0)) \\ \nonumber
& \leq   \norm[L^2]{\tilde{f}(0)- f_\infty}  +C  \left(\E_\lambda(\tilde{f}(0))-\E_\lambda(f_\infty)\right)^\theta \\ \label{Vacanza}
& \leq C \norm[C^2]{\tilde{f}(0)- f_\infty} ^{\theta} \, ,
\end{align} 
since $\theta <1$ and for $\varepsilon$ small enough. Then by real interpolation (see \cite[(i) of Thm. at page 29 and (ii) of Thm. at page 5]{Tri92},  \cite[Thm.~6.4.5(3)]{Lofstrom}
and \cite[Thm.4.6.1(e), page 328]{Triebel78}) there exists some $0 < \beta < 1$ such that for all $t \in (0,T)$ 
\begin{align*}
\norm[H^4]{\tilde{f}(t)- f_\infty} &  \leq C \| \tilde{f}(t)- f_\infty\|_{C^{4}([0,1])}\leq C \norm[L^2]{\tilde{f}- f_\infty}^{\beta}  \|\tilde{f}- f_\infty\|_{C^{4,\gamma}([0,1])}^{1-\beta} \, .
\end{align*}
By \eqref{Vacanza}, \eqref{514b}, and  \eqref{N0est}  we find for $\varepsilon$ small enough and for all $t \in (0,T)$
\begin{align*}
\norm[H^4]{\tilde{f}(t)- f_\infty} & \leq C \norm[C^2]{\tilde{f}(0)- f_\infty}^{\beta \theta} \leq C(\varepsilon)< \frac12 \sigma \, ,
\end{align*}
which gives a contradiction to the maximality of $T$ as in \eqref{defT} if $T <T'$. Thus $T=T'$ (with \eqref{defT} holding up to $T'$).  But then if $T' <\infty$, since \eqref{514b} holds up to time $T'$, we can start the flow again. Hence it must be $T'=\infty$ and \eqref{defT} holds for $t \in (0, \infty)$.
Finally notice that the flow is not just eternal but also smooth (see \cite[Theorem~8.1]{Poppenberg}).
Using \eqref{derGest} we immediately infer that $\norm[L^2]{\partial_t\tilde{f}} \in L^1(\R_{+})$. Furthermore from $\| \tilde{f}(t) - \tilde{f}(t') \|_{L^{2}} \leq \int _{t'}^{t} \| \partial_{t} \tilde{f} \|_{L^{2}}$
 it follows that there exists
 \begin{align}\label{limitevero}
 \tilde{f}_{\infty}:=\lim_{t \to \infty} \tilde{f}(t, \cdot) \mbox{ in }L^2 \, .
 \end{align}  

It remains to establish the relation between $f$ and $\tilde{f}$. We will see that $\tilde{f}$ is a reparametrization of $f$. From \eqref{solnormal} it follows that $\tilde{f}$ satisfies the differential equation 
\[ \partial_t \tilde{f} + \xi \partial_x \tilde{f}= - \nabla_{L^{2}} \mathcal{E}_{\lambda} (\tilde{f}) \mbox{ on } (0,\infty) \times I \, ,\]
for some smooth function $\xi$.  The tangential term can be generated via diffeomorphism. In fact, by 
\cite[Theorem~9.48]{Lee} 
for $t >0$ there exist smooth diffeomorphisms $\Phi(t,\cdot): I \to I$ such that
 \begin{equation}\label{DiffLee}
\left\{ \begin{array}{ll}
\partial_t \Phi= \xi \circ \Phi & \mbox{ in } (0,\infty)\times I,\\
\Phi(0, \cdot)= \id_{I} \, .
\end{array}
\right. 
\end{equation}  
Then the function $\tilde{f}(t, \Phi(t,x))$ satisfies the equation
\begin{align*}
\partial_t [\tilde{f}(t, \Phi(t,x))] & = (\partial_t \tilde{f})(t, \Phi(t,x)) +   (\partial_x \tilde{f})(t, \Phi(t,x)) \partial_t \Phi(t,x) \\
& = - \nabla_{L^{2}}\mathcal{E}_{\lambda} (\tilde{f})(t, \Phi(t,x))   - \xi(t,\Phi(t,x)) (\partial_x \tilde{f})(t, \Phi(t,x))  \\
& \quad +   (\partial_x \tilde{f})(t, \Phi(t,x)) \partial_t \Phi(t,x)  \\
& =  - \nabla_{L^{2}}\mathcal{E}_{\lambda} (\tilde{f})(t, \Phi(t,x))   = -  \nabla_{L^{2}}\mathcal{E}_{\lambda} (\tilde{f}(t, \Phi(t,x)))
\end{align*}
and the initial condition  
\[\tilde{f}(0, \Phi(0,x))= \tilde{f}(0, x) = \hat{f}_0 \circ \Phi_0(x)=f(t_k,\Psi^{-1}(t_k, \Phi_0(x))) \, .\]
Due to the uniqueness of the solution of the elastic flow (recall also \cite{MaMa}, \cite{KSnotes}) 
we find $f(t_k+t) \equiv \tilde{f}(t, \Phi \circ (\Psi^{-1}(t_k) \circ \Phi_0)^{-1})$, 
or equivalently
\[\tilde{f}(t,\cdot) = f(t_k+t,\Psi^{-1}(t_k)\circ \Phi_0 \circ \Phi^{-1}(t) (\cdot))\]
and, by \eqref{limitevero}, the claim follows.
\end{proof}
\begin{remark}\label{rem5.3}
Thanks to the  \L ojasiewicz-Simon gradient inequality we have also information about the rate of convergence. 
With the same notation as in the proof of Theorem \ref{main2} and computing as in \eqref{eq:ineqG} (recall also \eqref{6.8}) we find 

\begin{align*}
 -\frac{\diff}{\diff t} G(t) & = \theta G(t)^{\frac{(\theta-1)}{\theta}}\norm[L^2 (\diff s_{\tilde{f}})]{\nabla_{L^2} \E_\lambda(\tilde{f}(t))}^2 \\
  & \geq \theta G(t)^{\frac{(\theta-1)}{\theta}}\frac{1}{C^2} \left(\E_\lambda(\tilde{f}(t)) - \E_\lambda(f_\infty)\right)^{2(1-\theta)} 
\nonumber  = \frac{\theta}{C^2} G(t)^{\frac{(1-\theta)}{\theta}}, 
\end{align*}

from which it follows that $G(t)$ is $O(e^{-ct})$ for $\theta=\frac12$ and $O(t^{-\frac{\theta}{1-2\theta}})$  for the other values of $\theta$. By \eqref{derGest} we have
\begin{equation}\label{decadenza}
\norm[L^2]{\tilde{f}(t) - \tilde{f}_\infty} \leq \int_t^\infty \norm[L^2]{\partial_t \tilde{f}} \diff \tau = \frac{2C}{\theta} G(t) \, ,
\end{equation}
from which we derive the following rate of convergence
\begin{equation}
 \norm[L^2]{\tilde{f}(t) - \tilde{f}_{\infty}} \in 
      \left\{
      \begin{array}{ll}
                  O(e^{-ct})&\;\; \theta = \frac{1}{2}\\
                  O(t^{-\frac{\theta}{1-2\theta}})& \;\;\theta \in (0,\frac{1}{2})                 
      \end{array} \right.
\end{equation}
as $t \to \infty$.
\end{remark}

\appendix
\section{Corollary 3.11 in \texorpdfstring{\cite{Chill}}{Chi03}}\label{Chill} 
Our proof of the  \L ojasiewicz-Simon Inequality is based on \cite[Corollary 3.11]{Chill}.  For completeness we recall the setting in \cite{Chill} and then show how it is sufficient to  prove \textit{1.,2.,3.} in Remark \ref{rem:list} to apply \cite[Corollary 3.11]{Chill}.

In \cite{Chill} the following framework is considered. Let $V$ be a Banach space, $U \subset V$ open and $E\in C^2(U,\R)$ be the considered energy functional. Then the first order derivative of $E$ (denoted by $\mathcal{M}$) is in $C^1(U,V^{*})$ and the second order derivative of $E$ (denoted by $\mathcal{L}$) is in $C(U,B(V,V^{*}))$.  As usual $V^*$ denotes the dual of a Banach space $V$ and $B(X,Y)$ denote the Banach space of linear bounded operators from a Banach space $X$ to a Banach space $Y$. 

Let $\varphi$ be a critical point for $E$ in $U$, that is $\mathcal{M} (\varphi)=0$.\smallskip

\begin{hyp} (\cite[Hypothesis 3.2]{Chill})\label{hyp1}
  There exists a projection $P \in B(V,V)$ such that 
$$\Ima P = \Ker \mathcal{L}(\varphi)=: V_0 \, .$$
\end{hyp}

In this case we have $V=V_0 \oplus V_1$ (topological sum) with $V_{1}= \Ker P$. Let $P^* \in B(V^{*}, V^{*})$ be the adjoint projection. Then 
$$V^{*} = \Ima P^{*} \oplus  \Ker P^{*}=:   V_0^{*} \oplus   V_1^{*} \, .$$
Notice that this notation makes sense since we may identify $V_0^*$ with $\Ima P^*$ and $V_1^*$ with $\Ker P^*$.

\begin{hyp} (\cite[Hypothesis 3.4]{Chill}) \label{hyp2}
 There exists a Banach space $W$ such that
\begin{enumerate}
\item[(i)] $W \hookrightarrow V^{*}$ with continuous embedding;
\item[(ii)] Let $P$ be the projection from Hyp.~\ref{hyp1}. Then the adjoint $P^{*} \in B(V^{*},V^{*})$ leaves $W$ invariant; 
\item[(iii)] $\mathcal{M}\in C^1(U,W)$;
\item[(iv)] $\Ima\mathcal{L}(\varphi) = \Ker P^{*} \cap W$.
\end{enumerate}
\end{hyp}
We may state now Corollary 3.11 in \cite{Chill}.

\begin{corol}(\cite[Corollary 3.11]{Chill})
Let $\varphi$ be a critical point for $E$ in $U$ and assume Hypothesis~\ref{hyp1} and \ref{hyp2}. Assume in addition that there exist Banach spaces $X\subset V$ and $Y \subset W$ such that
\begin{enumerate}
\item the spaces $X$ and $Y$ are invariant under the projections $P$ and $P^{*}$ respectively;
\item the restriction of the derivative $\mathcal{M}$ to $U \cap X$ is analytic in a neighbourhood of $\varphi$ with values in $Y$;
\item $ \Ker \mathcal{L}(\varphi)$ is contained in $X$ and finite dimensional, and
\item $\Ima \left. \mathcal{L}(\varphi)\right|_{X}= \Ker P^{*} \cap Y$.
\end{enumerate}
Then the functional $E$ satisfies the \L ojasiewicz-Simon Inequality near $\varphi$,  that is there exist $\sigma>0$, $\theta 
\in (0, \frac{1}{2}]$, and $C \geq0$ such that for every $v \in U$ with $\|v -\varphi \|_{V} \leq \sigma $ 
$$|E(v)-E(\varphi)|^{1-\theta} \leq C\| \mathcal{M}(v) \|_{W}.$$ 
\end{corol}

We discuss now how, in our situation, it is sufficient to  prove \textit{1.,2.,3.} in Remark~\ref{rem:list} to apply \cite[Corollary 3.11]{Chill}.

We take: $X=V:= H^{4,\bot}_\text{c}$, $Y=W=H:= L^{2,\bot} $, $E$ as defined in Definition~\ref{def2} and   $\varphi=0$ as a critical point. That is, we consider the case that $\bar{f}$ is a critical point for the elastic energy $\mathcal{E}$.
To verify the assumptions of \cite[Corollary 3.11]{Chill} we take advantage of the Hilbert structure of $H$ as follows.

Since $V = H^{4,\bot}_\text{c}$ embeds densely into the Hilbert space $L^{2,\bot} = H$, we have the usual Gel\rq{}fand triple
\[
V \hookrightarrow H \cong H^* \hookrightarrow V^*.
\]
 
As discussed below Definition \ref{def2} the Fr\'echet-derivative $E\rq{}$ may be identified  with the $L^{2,\bot}$ - gradient of $E$, $d E: U \to H$ defined by 
\begin{equation} \label{eq:gradientE}
 E\rq{}(\psi)\phi = \langle d E(\psi), \phi \rangle _{L^{2}(\diff x)} \; \forall \psi \in U, \forall \phi \in V. \end{equation}

By \textit{2.} in Remark \ref{rem:list} $d E $ is analytic and we may consider its derivative
\[
  (d E)\rq{}: U \to \mathcal{B}(V,H).
\]

Let $\mathfrak{L} = (d E)\rq{}(0) : V \to H$. For $\phi$ and $\psi \in V$ we can write 
\begin{align}\nonumber
\langle \mathfrak{L}\phi,\psi\rangle_{L^2(dx)} 
&= \langle (d E)\rq{}(0) \phi, \psi \rangle_{L^2}
= \langle \frac{\diff}{\diff t} d E (t\phi) |_{t = 0}, \psi\rangle_{L^2}\\ \nonumber
&= \frac {\diff} {\diff t} \langle d E(t\phi), \psi \rangle_{L^2} |_{t=0}
= \frac{\diff}{\diff t} E\rq{}(t\phi)\psi|_{t = 0} \\ \label{ciccio}
&= E\rq{}\rq{}(0)(\phi,\psi)
\end{align}
which is symmetric in $\phi$ and $\psi$ since $E$ is analytic by \textit{1.} in Remark \ref{rem:list} and hence in particular $E \in C^2(U, \R)$ (an expression for the second variation of $E$ is provided  in Proposition~\ref{prop:2ndvar}). 

Set $$V_{0}:= \Ker \mathfrak{L} \subset V \subset H.$$  First of all consider $V_{0}$ as a subset of $H$, thus endowed with the $L^{2}$-scalar product.
If $\phi \in \Ker \mathfrak{L}$ and $\psi \in V$, then from $0= \langle \mathfrak{L} \phi,\psi \rangle_{L^{2}} = \langle \phi,  \mathfrak{L} \psi \rangle_{L^{2}} $ we infer that $\Ker \mathfrak{L}$ and $\Ima \mathfrak{L}$ are orthogonal subspaces in $H$, namely
\begin{equation}\label{ortogonality}
\Ker \mathfrak{L} \bot_{L^2} \Ima\mathfrak{L}. 
\end{equation}

Now we use 3. in Remark \ref{rem:list}, namely that $\mathfrak{L}:V \to H$ is Fredholm, that is $\text{dim} \Ker \mathfrak{L} < \infty$, $\Ima\mathfrak{L}$ is closed in $H$ and $\text{codim} \Ima\mathfrak{L} < \infty$ (see for instance \cite[page 168]{Brezis}). Then, $V_0$ (the kernel of $ \mathfrak{L}$) is finite dimensional and hence closed. There exists then a complementing subspace $V_1:= V_0^{\perp}$ such that 
$$ V  = V_{0} \oplus V_{1},$$
and a linear continuous projection $P:V \to V$ with $\Ima P=V_0$ (\cite[III.1 Thm.~1, Thm.~2]{Yosida}).
Thus Hypothesis~\ref{hyp1} (namely \cite[Hypothesis 3.2]{Chill}) is satisfied. Moreover 
 any $v \in V$ can be written as  $v =\tilde{v}_{0} + w$, with $\tilde{v}_{0} \in V_{0}$, $w \in \Ker P$.

 The adjoint operator $P^{*} \in \mathcal{B}(V^{*}, V^*)$ is again a linear and continuous projection and hence we can write 
  $$ V^{*}= \Ima P^{*} \oplus \Ker P^{*}= V_{0}^{*} \oplus V_{1}^{*}.$$
   We show now that $\Ima \mathfrak{L} = H^* \cap V_{1}^{*}$: indeed (after the canonical identification of $y^{*} \in H^{* } $  with $\langle y, \cdot \rangle_{L^{2}}$ for $y \in H$) we can write using the density of $V$ in $H$, the continuity of $P$, and~\eqref{ortogonality}
   \begin{align*}
   H^*\cap V_{1}^{*}& =\{  y^{*} \in H^{*} \cap V^*\, : \, P^{*}y^{*}(v)= y^{*}(Pv)=0 \quad \forall \, v \in V\} \\
  & = \{  y \in H \, : \, \langle y , P v \rangle_{L^{2}} =0 \quad  \forall \, v \in V \}\\
  & = \{  y \in H \, : \, \langle y , P x \rangle_{L^{2}} =0 \quad  \forall \, x \in H \} = (\Ima P)^{\perp}= V_{0}^{\perp}= \Ima \mathfrak{L} .
   \end{align*}
 
This gives   Hypothesis~\ref{hyp2} (iv)   
 and 4. in \cite[Cor.3.11]{Chill}. 
Next notice that
\begin{align*}
   H^*\cap V_{0}^{*}&=  \Ima P^{*} \cap H^* \\
&   =\{  y^{*} \in H^{*}\cap V^*\, : \, y^{*}=P^{*}v^{*}\quad \text{ for some } v^{*} \in V^{*}\} \\
& = \{ y^{*} \in H^{*}\cap V^* \, : \, y^{*}(x)= v^{*}(Px) \text{ for some } v^{*} \in V^{*} \text{ and } \, \forall \, x \in V \}.
   \end{align*}
Thus we infer   that  for  $y^{*} \in H^*\cap V_{0}^{*}$  and $w \in \Ker P \subset V$ we have  $y^{*}(w)=0$.

Finally we can show that $P^{*}$ leaves $H^*$ invariant in the sense that for $y^{*} \in H^{*} \cap V^{*}_0 $  (which we canonically identify with $\langle y, \cdot \rangle_{L^{2}}$ for $y \in H$) we have that $P^* y^*=y^*$. Indeed for  $v \in V$,
$v =\tilde{v}_{0} + w$, with $\tilde{v}_{0} \in V_{0}$, $w \in V_1= \Ker P$, we can write
$$ P^{*} y^{*}(v)=y^{*}(P v)= \langle y,P \tilde{v}_{0} \rangle_{L^{2}}=  \langle y,\tilde{v}_{0} \rangle_{L^{2}}= y^{*}(\tilde{v}_{0}) = y^{*}(\tilde{v}_{0})+ y^{*}(w) = y^{*}(v).$$
 We have just verified   Hypothesis~\ref{hyp2} (ii). 
 Assumptions   Hypothesis~\ref{hyp2} (i), (iii) 
 follow from the choice of the spaces and \textit{2.} in Remark \ref{rem:list}.

The assumptions of \cite[Corollary 3.11]{Chill} are also satisfied. Thus, by virtue of this corollary,  $E$ satifies the \L ojasiewicz-Simon inequality near $0$, i.e. there exists $\tilde{\sigma} \in (0,\rho)$, $\theta \in (0,\frac{1}{2}]$ and $C_2 > 0$ such that for all $\phi \in V, \norm[V]{\phi} < \tilde\sigma$,
\begin{equation} \label{eq:L.S}
 \betrag{E(0)-E(\phi)}^{1-\theta} \leq C_2 \norm[L^2(dx)]{d E(\phi)}.
\end{equation}

\section{Technical Proofs}  

\subsection{Analyticity}\label{appanaly}
\begin{proof}[Proof of Lemma \ref{lemma:analycitybyparts}]
\begin{enumerate}
  \item The map $U \to H^3(I, \R^{d})$, $\phi \mapsto \partial_x (\bar{f} + \phi)$ is well-defined, affine 
and continuous, hence it is analytic. Moreover, by definition of $U$ and since $H^3(I,\R^d)$ embeds into $C^0(I,\R^d)$, its image is contained in 
$$ V:= \{ \psi \in H^3(I, \R^d): \psi(x)\ne 0 \mbox{ for all }x \in [0,1]\} \, ,$$
that is an open set of $H^{3}(I, \R^{d})$. Since the Euclidean norm $| \cdot |: V \to H^{3}(I, \R)$ is analytic, we see that $F_1$ is given as the composition of two analytic mappings and hence analytic itself.
  \item Let consider the following open subset of $H^{3}(I, \R)$
$$ W:= \{ \psi \in H^3(I, \R): \psi(x)\ne 0 \mbox{ for all }x \in [0,1]\} \, .$$
Since $H^3(I,\R)$ embeds in $C^2([0,1],\R)$, for all $\psi \in W$ there exists $\delta_1=\delta_1(\psi)>0$ such that $|\psi(x)|\geq \delta_1>0$ for all $x \in [0,1]$ and hence $1/\psi$ is twice differentiable und its second derivative is given by
$$\left(\frac{1}{\psi}\right)'' = -\frac{\psi''}{\psi^2}+ 2\frac{(\psi')^2}{\psi^3} \, .$$
It follows from the definition of weak derivative that $1/\psi \in H^3(I,\R)$ and hence that the map
$$ G: W \rightarrow H^3(I,\R), \quad \psi \mapsto \frac{1}{\psi} \, ,$$
is well defined. We claim now that $G$ is also analytic. Let $\psi_0 \in W$ and $|\psi_0(x)|\geq \delta_1>0$ for all $x \in [0,1]$. We may assume w.l.o.g. that $\psi_0>0$ on $I$. By continuity, there is a $\delta_2 >0$ such that $\psi(x)\geq \frac12 \delta_1>0$ and $\|\psi-\psi_0\|_{\infty}\leq \frac12 \delta_1$ for all $\psi \in B_{\delta_{2}}(\psi_0) \subset H^3$. For such $\psi$ we write
$$ \frac{1}{\psi(x)}= \frac{1}{\psi_0(x)} \frac{1}{1-(1-\frac{\psi(x)}{\psi_0(x)})} = \frac{1}{\psi_0(x)} \sum_{k=0}^{\infty} \left(1-\frac{\psi(x)}{\psi_0(x)} \right)^k \, . $$
Since there is a universal constant $c$ such that $\| g_1 g_2 \|_{H^3} \leq c \| g_1 \|_{H^3} \| g_2 \|_{H^3}$, for all $g_1,g_2 \in H^3$, we can write,
$$G(\psi)=  \sum_{k=0}^{\infty}\frac{1}{(\psi_0)^{k+1}} \left(\psi_0- \psi \right)^k$$ for $\psi$ such that 
$$\| \psi- \psi_0\|_{H^3} < c^{-1}\| \frac{1}{\psi_0}\|^{-1}_{H^3} \mbox{ and }\|\psi-\psi_0\|_{H^3}< \delta_2\, .$$
This gives first the analyticity of $G$ in $\psi_0$ and then the analyticity of $G$ in $W$ since $\psi_0$ was arbitrary. \\
Since the map $F_1$ has values in $W$, $G \circ F_1:U \to H^3(I,\R)$ is analytic. Since $H^3(I,\R)$ is a Banach algebra and $F_2(\phi)= G(F_1(\phi)) (\partial_x (\bar{f} + \phi))$, we see that $F_2$ is an analytic function as product of analytic functions. 
  \item The map from $U$ to $H^2(I,\R^d)$ that associates to $\phi$ the vector field $\partial_x F_2(\phi)$ is analytic as composition of analytic mappings. With the same notation as in the previous part of the proof, since $H^3$ embeds continuously into $H^2$ also the function $G \circ F_1: U \to H^2(I,\R)$ is analytic. 
Since $H^2(I,\R)$ is a Banach algebra, it follows that the mapping
$$F_3(\phi) = G(F_1(\phi)) \partial_x F_2 (\phi) \, $$
is well defined and analytic.
  \item Due to the continuity of the embedding $H^2 \hookrightarrow L^2$, the mapping  $\phi \mapsto F_3(\phi)$, from $U$ to $L^2(I,\R^d)$ is analytic. Since $|\cdot|^2: L^2(I;\R^d) \to L^1(I,\R)$ is also analytic, then the mapping 
$$U \ni \phi \mapsto |F_3(\phi)|^2 \in  L^1(I,\R)\, ,$$
is analytic as composition of analytic functions. The mapping $U \ni \phi \mapsto F_1(\phi) \in C^0(I,\R)$ is also analytic. Using that the product of a continuous function with a $L^1$-function, is a $L^1$-function and that this product is bilinear and continuous we find that 
$$F_4: U \to L^1(I,\R), \quad \phi \mapsto |F_3(\phi)|^2 F_1(\phi) \, ,$$
is analytic. 
\item Since the Euclidean scalar product induces a bilinear continuous product from $H^2(I,\R^d) \times H^2(I,\R^d)$ into $H^2(I,\R)$ one sees that the mapping  $\phi \mapsto |F_3(\phi)|^2$, from $U$ to $H^2(I,\R)$ is analytic. Reasoning similarly and since $H^3 \hookrightarrow H^2$, one gets the analyticity of the map $\phi \mapsto |F_3(\phi)|^2 F_3(\phi) F_1(\phi)$ from $U$ to $H^2(I,\R^d)$. Due to the embedding  $H^2 \hookrightarrow L^2$, it follows that  $\frac{1}{2} |\Kapp|^2 \Kapp |\bar{f}_x + \phi_x|$ is analytic from $U$ to $L^{2}$.

To show that $ |\bar{f}_x + \phi_x|\nabla_s^2 \Kapp: U \to L^{2,\bot}$ is analytic, we write it explicitely in order to see it as product of analytic functions. We have
\begin{align*}
\partial_s \Kapp & = \frac{1}{|\partial_x (\bar{f}+\phi)|} \partial_x \Kapp \, , \\
\nabla_s \Kapp & = \frac{1}{|\partial_x (\bar{f}+\phi)|} \partial_x \Kapp + | \Kapp|^2 \partial_s (\bar{f}+\phi)\, , \\
\partial_s \nabla_s \Kapp & = \frac{1}{|\partial_x (\bar{f}+\phi)|^2} \partial_x^2 \Kapp - \frac{1}{|\partial_x (\bar{f}+\phi)|^4}  \langle \partial_x (\bar{f}+\phi), \partial_x^2 (\bar{f}+\phi)\rangle_{\text{euc}} \partial_x \Kapp \\
& \quad + | \Kapp|^2 \Kapp + 2 \langle \Kapp, \partial_s \Kapp\rangle_{\text{euc}} \partial_s (\bar{f}+\phi)\, 
\end{align*}
from which it follows that 
\begin{align*}
& \qquad  |\bar{f}_x + \phi_x|\nabla_s^2 \Kapp  \\
& = \frac{1}{|\partial_x (\bar{f}+\phi)|} \partial_x^2 \Kapp - \frac{1}{|\partial_x (\bar{f}+\phi)|} \langle \partial_x^2 \Kapp, \partial_s (\bar{f}+\phi)\rangle_{\text{euc}}  \partial_s (\bar{f}+\phi) \\
& \quad  + | \Kapp|^2 \Kapp  |\bar{f}_x + \phi_x|  - \frac{1}{|\partial_x (\bar{f}+\phi)|^3}  \langle \partial_x (\bar{f}+\phi), \partial_x^2 (\bar{f}+\phi)\rangle_{\text{euc}} \partial_x \Kapp \\
& \quad  - \frac{1}{|\partial_x (\bar{f}+\phi)|^2}  \langle \partial_x (\bar{f}+\phi), \partial_x^2 (\bar{f}+\phi)\rangle_{\text{euc}} | \Kapp|^2  \partial_s (\bar{f}+\phi) \, .
\end{align*}
The analyticity of this mapping from $U$ to $L^2$ is established with the arguments used in the previous claims. Hence $F_5$ is analytic as sum of analytic functions.\qedhere
\end{enumerate}
\end{proof}

\subsection{Calculation of the second variation}
\begin{proof}[Proof of Proposition \ref{prop:2ndvar}] \label{proof:2ndvar} 
For simplicity of notation in this proof we denote the Euclidean scalar-product simply by $\metric{\cdot,\cdot}$ and we write $s=s_{\bar f}$ and $\Kapp=\Kapp_{\bar{f}}$. For $\phi,\psi \in U \subset H^{4,\bot}_c$ we find using \eqref{gradfat}
\begin{align}
E''(0) [\phi,\psi] & = \left. \frac{d}{d\varepsilon} \right|_{\varepsilon=0} \left. \frac{d}{d\sigma} \right|_{\sigma=0}\mathcal{E} (\bar{f}+ \varepsilon \phi + \sigma \psi) \nonumber \\
& =  \left. \frac{d}{d\varepsilon} \right|_{\varepsilon=0} \left. \frac{d}{d\sigma} \right|_{\sigma=0}\mathbf{E} ( \varepsilon \phi + \sigma \psi) 
 =  \left. \frac{d}{d\varepsilon} \right|_{\varepsilon=0} \mathbf{E}' ( \varepsilon \phi )( \psi) \nonumber \\
& = \left. \frac{d}{d\varepsilon} \right|_{\varepsilon=0} \langle \nabla_{s_{\bar{f} + \varepsilon \phi }}^2 \Kapp_{\bar{f} + \varepsilon \phi } + \frac{1}{2} |\Kapp_{\bar{f} + \varepsilon \phi }|^2 \Kapp_{\bar{f} + \varepsilon \phi }, \psi\rangle_{L^2(\diff s_{\bar{f}+ \varepsilon \phi })} \, . \label{formulagiusta}
\end{align}
We write $s_{\varepsilon}:=s_{\bar{f} + \eps \phi}$ and $\Kapp_{\varepsilon}: =\Kapp_{\bar{f} + \varepsilon \phi }$. Since $\partial_{s} = |\partial_x \bar{f}|^{-1}\partial_x$, one finds that
\begin{align}
\deps |\partial_x(\bar{f}+\eps \phi)| \Big|_{\eps = 0}&=  \deps \metric{ \partial_x (\bar{f}+\eps \phi), \partial_x(\bar{f} + \eps \phi)}^\frac{1}{2} \Big|_{\eps = 0}\nonumber\\
 &= \frac{1}{2 |\partial_x \bar{f}|} 2 \metric{\partial_x  \phi, \partial_x \bar{f}} = \metric{\partial_x  \phi,\diffs \bar{f}} =  -|\partial_x \bar{f}| \metric{\phi,\Kapp}. \label{eq:varnorm}
\end{align}
From this and the fact that $\metric{\phi,\diffs \bar{f}} = 0$ is constant we see that the variation of the volume form $\diff s$ is given by
\begin{align}
\deps{\diff s_\eps } \Big|_{\eps = 0}&=  - \metric{\phi,\Kapp} \diff s.\label{eq:varform}
\end{align}
Using \eqref{eq:varnorm} again one also finds that the variation of the derivative with respect to arc length is given by
\begin{align}
\deps \left(\diffseps\eta\right)   \Big|_{\eps = 0} &= \frac{\partial_x(\deps {\eta})}{|\partial_x \bar{f}|} \Big|_{\eps = 0} +( \partial_x \eta)  |_{\eps = 0}\cdot \left(- \frac{1}{|\partial_x \bar{f}|^2}|\partial_x \bar{f}|(-\metric{\phi,\Kapp })\right)\nonumber\\
&{  = (\partial_s \partial_{\varepsilon} \eta) |_{\eps=0} +\metric{ \phi, \Kapp}  (\diffs \eta)|_{\eps=0} },  \label{eq:vardiffseps}
\end{align}
where $\eta$ is any sufficiently smooth function on $I \times (-\delta, \delta)$. This allows us to calculate
\begin{align}
\deps \diffseps(\bar{f}+\eps \phi)  \Big|_{\eps = 0}&= \diffs\phi + \metric{\phi,\Kapp}\diffs \bar{f} { = \nabla_s \phi} \label{eq:vardiffsepsphi}
\end{align}
and
\begin{align}
\deps{\Kapp_\eps}  \Big|_{\eps = 0}&=\deps  \left( \diffseps^2 (\bar{f}+\eps \phi) \right) \Big|_{\eps = 0} \nonumber \\
&\stackrel{\eqref{eq:vardiffseps}}{=} \diffs\left[(\diffseps(\bar{f}+\eps \phi)) \right]  \Big|_{\eps = 0} + \metric{\phi,\Kapp}\diffs (\diffs \bar{f}) \nonumber \\
&\stackrel{\eqref{eq:vardiffsepsphi}}{=} \diffs\nabla_s \phi  -\metric{\diffs \nabla_s \phi, \diffs \bar{f}}\diffs \bar{f} +\metric{\diffs \nabla_s \phi, \diffs \bar{f}}\diffs \bar{f}+ \metric{\phi,\Kapp}\Kapp \nonumber\\
&= \nabla_s^2 \phi- \metric{\nabla_s \phi, \Kapp} \diffs \bar{f} + \metric{\phi, \Kapp} \label{eq:varkapp} \Kapp 
\end{align}
where we used the definition of $\nabla_s$ and the orthogonality in the last line. Hence for $\psi \in U$ we get
\allowdisplaybreaks{\begin{align}
 & \qquad \left. \frac{d}{d\varepsilon} \right|_{\varepsilon=0}  \left( \frac{1}{2} \int_I |\Kapp_\eps|^2 \metric{\Kapp_\eps, \psi} \diff s_\eps    \right ) \nonumber \\
  &\stackrel{\eqref{eq:varform}}{=} \frac{1}{2} \int_I 2\metric{\deps{\Kapp_\eps}, \Kapp}  \Big|_{\eps = 0} \metric{\Kapp,\psi}+ |\Kapp|^2 \metric{\deps{\Kapp_\eps}, \psi}  \Big|_{\eps = 0} -|\Kapp|^2 \metric{\Kapp, \psi} \metric{\Kapp,\phi} \diff s \nonumber \\
 &\stackrel{\eqref{eq:varkapp}}{=}\int_I \left(\metric{ \nabla_s^2 \phi,\Kapp} +\metric{\metric{\phi, \Kapp} \Kapp, \Kapp}\right)\metric{\Kapp,\psi}  \nonumber \\
&\qquad  \qquad + \frac{1}{2} |\Kapp|^2 (\metric{\nabla_s^2 \phi,\psi}+ \metric{\phi, \Kapp}\metric{ \Kapp, \psi})  -\frac{1}{2}|\Kapp|^2 \metric{\Kapp, \psi} \metric{\Kapp,\phi}  \diff s \nonumber \\
&=   \int_I \metric{ \nabla_s^2 \phi,\Kapp}\metric{\Kapp,\psi} +|\Kapp|^2 \metric{\Kapp, \psi} \metric{\Kapp,\phi}   + \frac{1}{2} |\Kapp|^2 \metric{\nabla_s^2 \phi,\psi}\diff s \nonumber\\
&=   \int_I \metric{ \nabla_s^2 \phi,\Kapp}\metric{\Kapp,\psi} +|\Kapp|^2 \metric{\Kapp, \psi} \metric{\Kapp,\phi}   \nonumber\\
& \quad -\metric{\nabla_s \Kapp, \Kapp} \metric{\nabla_s \phi, \psi}
 - \frac{|\Kapp|^2}{2}\metric{\nabla_s\phi, \nabla_s \psi} \diff s \label{eq:varrechtInt}
 \end{align}}
where the last equality stems from the following calculation:
\begin{align}
 &\mkern-24mu  \int_I\frac{1}{2} |\Kapp|^2 \metric{\nabla_s^2 \phi,\psi} \diff  s \nonumber \\
 &=   \int_I\frac{1}{2} |\Kapp|^2 \metric{\partial_s\nabla_s \phi,\psi}\diff s \nonumber\\
  &=   \int_I\frac{1}{2} |\Kapp|^2 \left[ \left(\diffs \metric{\nabla_s \phi,\psi}\right) - \metric{\nabla_s \phi, \nabla_s \psi} \right] \diff s\nonumber\\
  &= 
  \left. \frac{1}{2} |\Kapp|^2 \metric{\nabla_s \phi,\psi}\right|_{\partial I}
  -\int_I  \metric{\diffs \Kapp, \Kapp} \metric{\nabla_s \phi, \psi}+ \frac{1}{2} |\Kapp|^2 \metric{\nabla_s \phi, \nabla_s \psi} \diff s\nonumber\\
  &= -\int_I \metric{\nabla_s \Kapp, \Kapp} \metric{\nabla_s \phi, \psi} + \frac12 |\Kapp|^2 \metric{\nabla_s\phi, \nabla_s \psi} \diff s \nonumber \,,
 \end{align}
since the boundary term vanishes in the last line as $\psi$ was assumed to be in $H^1_0$. Furthermore by definition of $\nabla_s$ and application of equations \eqref{eq:vardiffseps} and \eqref{eq:vardiffsepsphi} in the second line one finds
\begin{align}
\deps \left(\nabla_{s_\eps} \eta_\eps \right)  \Big|_{\eps = 0} &= \deps \diffseps \eta_\eps   \Big|_{\eps = 0} - \deps  \metric{\diffseps \eta_\eps,\diffseps (\bar{f} +\eta_\eps \phi)}  \diffseps (\bar{f} + \eps \phi) \Big|_{\eps = 0}  \nonumber\\
        &= \nabla_s \deps \eta_\eps \Big|_{\eps = 0}   
	   + \metric{\phi, \Kapp} \diffs\eta_\eps \Big|_{\eps = 0}
-\metric{\metric{\phi, \Kapp} \diffs \eta_\eps, \diffs \bar{f}} \diffs \bar{f} \Big|_{\eps = 0} \nonumber\\
	   &\qquad -\metric{\diffs \eta_\eps, \nabla_s \phi}\diffs \bar{f}\Big|_{\eps = 0}  -\metric{\diffs \eta_\eps, \diffs \bar{f}}\nabla_s \phi \Big|_{\eps = 0} \nonumber\\
	     &=\left( \nabla_s \left(\deps{\eta_\eps}\right)   
	   + \metric{\phi, \Kapp} \nabla_s \eta_\eps  -\metric{\nabla_s \eta_\eps, \nabla_s \phi}\diffs \bar{f}  + \metric{\eta_\eps, \Kapp} \nabla_s \phi \right) \Big|_{\eps = 0} \label{eq:varnabla}
\end{align}
for {$\eta_\eps$ which is normal to $\bar{f}$ at $\eps = 0$}. Hence one finds that for $\psi$ which is not depending on $\eps$ and orthogonal to $\diffs \bar{f}$, putting $\eta_\eps := \nabla_{s_\eps} \psi$ in the first step and $\eta = \psi$ in the second:
\begin{align}
\deps \nabla_{s_\eps}^2 \psi \Big|_{\eps = 0}   &= \deps  \nabla_{s_\eps}  \left(\nabla_{s_\eps}\psi\right) \Big|_{\eps = 0}  \nonumber\\
&\stackrel{\eqref{eq:varnabla}}{=} \nabla_s [\deps \nabla_{s_\eps} \psi] \Big|_{\eps = 0} + \metric{\phi, \Kapp} \nabla_s \left(\nabla_{s}\psi\right)\nonumber \\
    & \qquad \qquad -\metric{\nabla_s \nabla_s\psi, \nabla_s \phi}\diffs \bar{f} + \metric{\nabla_s\psi,\Kapp} \nabla_s \phi\nonumber\\\begin{split}
 & \stackrel{\eqref{eq:varnabla}}{=} \nabla_s
      \left[\metric{\phi, \Kapp} \nabla_s \psi -\metric{\nabla_s\psi, \nabla_s \phi}\diffs \bar{f} + \metric{\psi, \Kapp} \nabla_s \phi
    \right]   \\ 
     &\qquad \qquad+ \metric{\phi, \Kapp} \nabla_s^2\psi -\metric{\nabla_s^2 \psi, \nabla_s \phi}\diffs \bar{f} + \metric{\nabla_s\psi,\Kapp} \nabla_s \phi. \label{eq:varnablasquared}
     \end{split}
\end{align}
To apply these formulas in the calculation of the second derivation, we need the following formulas of partial integration with respect to $\nabla$ for arbitrary (not necessarily orthogonal) functions $\eta$ and $\xi$ which are differentiable once or twice, respectively):\begin{align}
{  \int_I \metric{\nabla_s \eta, \xi}\diff s} &= \int_I \metric{\diffs \eta, \xi}- \metric{\diffs \bar{f}, \xi}\metric{\diffs \eta, \diffs \bar{f}} \diff s \nonumber\\
&= \int_I \left(\diffs \metric{\eta,\xi}\right) - \metric{\eta, \diffs \xi}- \metric{\diffs \bar{f}, \xi}\metric{\diffs \eta, \diffs \bar{f}} \diff s \nonumber\\
&=\left.\metric{\eta,\xi}\right|_{\partial I}+\int_I- \metric{\eta, \nabla_s \xi + \metric{\diffs \xi, \diffs \bar{f}} \diffs \bar{f}}- \metric{\diffs \bar{f}, \xi}\metric{\diffs \eta, \diffs \bar{f}} \diff s\nonumber\\
&{=\left.\metric{\eta,\xi}\right|_{\partial I}-\int_I \metric{\eta, \nabla_s \xi} + \metric{\diffs \xi, \diffs \bar{f}}\metric{\eta, \diffs \bar{f}}+ \metric{\diffs \bar{f}, \xi}\metric{\diffs \eta, \diffs \bar{f}} \diff s.} \label{eq:partintnabla}
\end{align}
And thus
\begin{align}
 {  \int_I \metric{\nabla_s^2 \eta, \xi}\diff s} &=  \int_I \metric{\nabla_s (\nabla_s\eta), \xi}\diff s\nonumber\\
  &\stackrel{\eqref{eq:partintnabla}}{=} \left.\metric{\nabla_s \eta,\xi}\right|_{\partial I}\nonumber\\
  &\qquad +\int_I- \metric{\nabla_s\eta, \nabla_s \xi} - \metric{\diffs \xi, \diffs \bar{f}}\cdot 0- \metric{\diffs \bar{f}, \xi}\metric{\diffs \nabla_s\eta, \diffs \bar{f}} \diff s.\nonumber\\
  &\stackrel{\eqref{eq:partintnabla}}{=}\left.\metric{\nabla_s \eta,\xi} - \metric{\eta,\nabla_s\xi}\right|_{\partial I} \nonumber \\
   &\qquad + \int_I -\left( -\metric{\eta, \nabla_s^2 \xi} - \metric{\diffs\nabla_s \xi, \diffs \bar{f}}\metric{\eta, \diffs \bar{f}}- \metric{\diffs \bar{f},\nabla_s \xi}\metric{\diffs \eta, \diffs \bar{f}} \right)\nonumber \\
  &\qquad  \qquad - \metric{\diffs \bar{f}, \xi}\metric{\diffs \nabla_s\eta, \diffs \bar{f}} \diff s\nonumber\\
  \begin{split}
     &  = \left.\metric{\nabla_s \eta,\xi} - \metric{\eta,\nabla_s\xi}\right|_{\partial I} \\
  &\qquad {  + \int_I\metric{\eta, \nabla_s^2 \xi} + \metric{\diffs\nabla_s \xi, \diffs \bar{f}}\metric{\eta, \diffs \bar{f}} - \metric{\diffs \bar{f}, \xi}\metric{\diffs \nabla_s\eta, \diffs \bar{f}} \diff s} \label{eq:partintnablasquared}
  \end{split}
\end{align}
Furthermore it holds for the normal vector fields $\psi, \phi$ by \eqref{eq:vardiffsepsphi}, \eqref{eq:varkapp} and orthogonality:
 \begin{align}
 &{ \deps \metric{\diffseps(\bar{f}+\eps \phi), \Kapp_\eps } \metric{\diffseps\nabla_{s_\eps} \psi, \diffseps(\bar{f}+\eps \phi)}} \Big|_{\varepsilon=0}\label{eq:varrest2} \\
 &\qquad = {\metric{\nabla_s\phi, \Kapp}\metric{\diffs\nabla_s \psi, \diffs \bar{f}}} -\metric{\diffs \bar{f}, \diffs \bar{f}} \metric{\nabla_s \phi, \Kapp}\metric{\diffs \nabla_s \psi, \diffs \bar{f}}   = 0, \nonumber\, ,\\
& {  \deps \metric{\diffseps(\bar{f}+\eps \phi), \psi }\metric{\diffseps\nabla_{s_\eps} \Kapp_\eps, \diffseps(\bar{f}+\eps \phi)}}   \Big|_{\varepsilon=0} \nonumber\\
&\qquad \stackrel{\eqref{eq:vardiffsepsphi}}{=} \metric{\nabla_s \phi, \psi}\metric{\partial_s \nabla_s \Kapp, \diffs \bar{f}} 
= {  -\metric{\nabla_s \phi, \psi}\metric{\nabla_s \Kapp, \Kapp}} \label{eq:varrest1}
 \end{align}
 and
 \begin{align}
  \deps{\metric{\Kapp_\eps, \nabla_{s_\eps}^2 \psi }}  \Big|_{\varepsilon=0} &= \metric{\deps{\Kapp_\eps}, \nabla_s^2 \psi}  \Big|_{\varepsilon=0}+ \metric{\Kapp,\deps{\nabla_{s_\eps} ^2 \psi}} \Big|_{\varepsilon=0}\nonumber\\ 
& \!\!\stackrel{\eqref{eq:varkapp}}{=} \metric{\nabla_s^2 \phi, \nabla_s^2 \psi} + \metric{\phi,\Kapp}\metric{\Kapp, \nabla_s^2 \psi}\nonumber\\
& \qquad \stackrel{\eqref{eq:varnablasquared}}{+} \metric{\Kapp,\nabla_s \left[ \metric{\phi, \Kapp} \nabla_s \psi -\metric{\nabla_s\psi, \nabla_s \phi}\diffs \bar{f} + \metric{\psi, \Kapp} \nabla_s \phi\right]}\nonumber\\
 & \,\,\,\,\qquad + \,\,\,\, \metric{\phi,\Kapp}\metric{\Kapp, \nabla_s^2 \psi} + \metric{\nabla_s \psi, \Kapp}\metric{\Kapp,\nabla_s \phi}. \label{eq:varlinkesintinnen1}
  \end{align}
Using \eqref{eq:partintnabla} we note putting $ \xi = \Kapp$ and $ \eta =\metric{\phi,\Kapp}\nabla_s \psi +\metric{\psi,\Kapp}\nabla_s \phi$ since $\eta$ vanishes on the boundary of $I$ that 
\begin{align}
 & \qquad \int_I\metric{\Kapp,\nabla_s \left[ \metric{\phi, \Kapp} \nabla_s \psi -\metric{\nabla_s\psi, \nabla_s \phi}\diffs \bar{f} + \metric{\psi, \Kapp} \nabla_s \phi\right]}  \diff s  = \nonumber \\
&=\int_I \metric{\Kapp,\nabla_s \left[ \metric{\phi, \Kapp} \nabla_s \psi + \metric{\psi, \Kapp} \nabla_s \phi\right]}  \diff s + \int_I \metric{\Kapp,\diffs \left[ -\metric{\nabla_s\psi, \nabla_s \phi}\diffs \bar{f} \right]}  \diff s\nonumber\\
& \stackrel{\eqref{eq:partintnabla}}{=}-\int_I \metric{\nabla_s\Kapp,\nabla_s\psi} \metric{\phi, \Kapp}  +\metric{\nabla_s\Kapp,\nabla_s\phi} \metric{\psi, \Kapp} + (-\underbrace{\metric{\diffs \Kapp, \diffs \bar{f}}}_{=-|\Kapp|^2}\metric{\nabla_s \psi, \nabla_s \phi})\diff s\nonumber
\end{align}
and hence by \eqref{eq:varlinkesintinnen1}
\begin{align}\nonumber
  \int_I\deps{\metric{\Kapp_\eps, \nabla_{s_\eps}^2 \psi }}  \Big|_{\varepsilon=0} \diff s 
& = \int_I \metric{\nabla_s^2 \phi, \nabla_s^2 \psi} + \metric{\phi,\Kapp}\metric{\Kapp, \nabla_s^2 \psi}  -\metric{\nabla_s\Kapp,\nabla_s\psi} \metric{\phi, \Kapp}\\ \nonumber
& \qquad   -\metric{\nabla_s\Kapp,\nabla_s\phi} \metric{\psi, \Kapp}    - |\Kapp|^2\metric{\nabla_s \psi, \nabla_s \phi}\\
& \qquad + \metric{\phi,\Kapp}\metric{\Kapp, \nabla_s^2 \psi} + \metric{\nabla_s \psi, \Kapp}\metric{\Kapp,\nabla_s \phi} \diff s \label{eq:varlinkesintinnen2}.
\end{align}

Now we use \eqref{eq:partintnablasquared} with $\eta := \Kapp_\eps, \xi := \psi, \bar{f} = \bar{f} + \eps \phi$. Due to the choice of the space the boundary terms disappear. Indeed, $\left.\metric{\nabla_{s_\eps} \Kapp_\eps, \psi}\right|_{\partial I}$ vanishes since $\psi \in H^1_0$ and for the derivative of the other boundary term we find
\begin{align}
\left. \partial_\eps \metric{ \Kapp_\eps, \nabla_{s_\eps} \psi}\right|_{\varepsilon = 0} &= \metric{\nabla_s^2 \phi - \metric{\nabla_s \phi, \Kapp} \partial_s \bar{f} + \metric{\phi, \Kapp} \Kapp, \nabla_s \psi}\nonumber \\
  &\;\;\; +  \metric{\Kapp,\metric{\phi,\Kapp} \nabla_s \psi - \metric{\nabla_s \psi,\nabla_s \phi}\partial_s \bar{f} + \metric{\psi,\Kapp}\nabla_s \phi}\nonumber\\
  &= \metric{\nabla_s^2 \phi + \metric{\phi,\Kapp}\Kapp, \nabla_s \psi} + \metric{\Kapp, \metric{\phi,\Kapp}\nabla_s \psi + \metric{\psi,\Kapp}\nabla_s \phi}\label{eq:boundarysecondvar}
\end{align}
by \eqref{eq:varkapp}, \eqref{eq:varnabla} and orthogonality. Since $\psi, \phi \in H^2_0$ 
 by definition of $H^{4,\bot}_{\text{c}}$ we find that \eqref{eq:boundarysecondvar} vanishes on the boundary of $I$. Hence
\begin{align}
& \left. \frac{d}{d\varepsilon} \right|_{\varepsilon=0} \int_I \metric{\nabla_{s_\eps}^2\Kapp_\eps,\psi } \diff s_\eps \nonumber\\
  & \stackrel{\eqref{eq:partintnablasquared}}{=}  \left. \frac{d}{d\varepsilon} \right|_{\varepsilon=0} \Big ( \int_I \metric{\Kapp_\eps, \nabla_{s_\eps}^2 \psi} + \metric{\diffseps(\bar{f}+\eps \phi), \Kapp_\eps}\metric{\diffseps\nabla_{s_\eps}\psi, \diffseps(\bar{f}+\eps \phi)} \nonumber\\
  & \qquad - \metric{\diffseps(\bar{f}+\eps \phi), \psi}\metric{\diffseps\nabla_{s_\eps}\Kapp_\eps,\diffseps(\bar{f}+\eps \phi)} \diff s_\eps \Big ) 
  \nonumber\\
  &\stackrel{\eqref{eq:varform}}{=} \int_I (\metric{\Kapp, \nabla_s^2 \psi}) (-\metric{\phi, \Kapp}) \diff s  + \int_I\deps{\metric{\Kapp_\eps, \nabla_{s_\eps}^2 \psi }}  \Big|_{\varepsilon=0}\diff s\nonumber\\
  &\qquad  + \int_I \deps{\metric{\diffseps(\bar{f}+\eps \phi), \Kapp_\eps}\metric{\diffseps\nabla_{s_\eps}\psi, \diffseps(\bar{f}+\eps \phi)}}  \Big|_{\varepsilon=0}\diff s  \nonumber\\
  &\qquad - \int_I \deps{\metric{\diffseps(\bar{f}+\eps \phi), \psi}\metric{\diffseps\nabla_{s_\eps}\Kapp_\eps,\diffseps(\bar{f}+\eps \phi)}}  \Big|_{\varepsilon=0} \diff s\nonumber\\
  &\stackrel{}{=} \int_I -\metric{\Kapp, \nabla_s^2 \psi}\metric{\phi, \Kapp} \diff s + \stackrel{\eqref{eq:varlinkesintinnen2}}{+} \int_I \metric{\nabla_s^2 \phi, \nabla_s^2 \psi} + \metric{\phi,\Kapp}\metric{\Kapp, \nabla_s^2 \psi}\nonumber\\
&\qquad  -\metric{\nabla_s\Kapp,\nabla_s\psi} \metric{\phi, \Kapp}  -\metric{\nabla_s\Kapp,\nabla_s\phi} \metric{\psi, \Kapp} - |\Kapp|^2\metric{\nabla_s \psi, \nabla_s \phi}\nonumber\\
&\qquad + \metric{\phi,\Kapp}\metric{\Kapp, \nabla_s^2 \psi} + \metric{\nabla_s \psi, \Kapp}\metric{\Kapp,\nabla_s \phi} \diff s \stackrel{\eqref{eq:varrest2}}{+} 0\nonumber \\
  &\quad \stackrel{\eqref{eq:varrest1}}{-}\int_I -\metric{\nabla_s \phi, \psi}\metric{\nabla_s \Kapp, \Kapp} \diff s\nonumber\\ \begin{split}                                                                                                                        
  &=  \int_I \metric{\nabla_s^2 \phi, \nabla_s^2 \psi} -\metric{\nabla_s\Kapp,\nabla_s\psi} \metric{\phi, \Kapp}  -\metric{\nabla_s\Kapp,\nabla_s\phi} \metric{\psi, \Kapp}\\
& \qquad \qquad   - |\Kapp|^2\metric{\nabla_s \psi, \nabla_s \phi}+ \metric{\phi,\Kapp}\metric{\Kapp, \nabla_s^2 \psi} + \metric{\nabla_s \psi, \Kapp}\metric{\Kapp,\nabla_s \phi}\\
& \qquad \qquad +\metric{\nabla_s \phi, \psi}\metric{\nabla_s \Kapp, \Kapp}\diff s \label{eq:varlinkesint}.\end{split}
\end{align}
From equations \eqref{formulagiusta}, \eqref{eq:varlinkesint} and \eqref{eq:varrechtInt} we finally have
\begin{align}
 E''(0)(\phi,\psi)   &= 
 \left. \frac{d}{d\varepsilon} \right|_{\varepsilon=0}  \Big(
 \int_I \metric[\text{euc}]{\nabla_{s_\eps}^2\Kapp_\eps, \psi} \diff s_{\eps} 
 + \frac{1}{2} 
 \int_I \metric{\Kapp_\eps,\Kapp_\eps}\metric{\Kapp_\eps,\psi} \diff s_{\eps}   \Big)
 \nonumber\\
&= \int_I \metric{\nabla_s^2 \phi, \nabla_s^2 \psi}- \metric{\nabla_s\Kapp,\nabla_s\psi} \metric{\phi, \Kapp}  -\metric{\nabla_s\Kapp,\nabla_s\phi} \metric{\psi, \Kapp}\nonumber\\
&\qquad  - |\Kapp|^2\metric{\nabla_s \psi, \nabla_s \phi}+ \metric{\phi,\Kapp}\metric{\Kapp, \nabla_s^2 \psi}\nonumber\\
&\qquad + \metric{\nabla_s \psi, \Kapp}\metric{\Kapp,\nabla_s \phi}+\metric{\nabla_s \phi, \psi}\metric{\nabla_s \Kapp, \Kapp}\diff s\nonumber\\
& \qquad + \int_I \metric{ \nabla_s^2 \phi,\Kapp}\metric{\Kapp,\psi} +|\Kapp|^2 \metric{\Kapp, \psi} \metric{\Kapp,\phi}\nonumber\\
&\qquad  -\metric{\nabla_s \Kapp, \Kapp} \metric{\nabla_s \phi, \psi} - \frac{|\Kapp|^2}{2}\metric{\nabla_s\phi, \nabla_s \psi} \diff s \nonumber\\
\begin{split} \nonumber
 &= \int_I \metric{\nabla_s^2 \phi, \nabla_s^2 \psi}- \metric{\nabla_s\Kapp,\nabla_s\psi} \metric{\phi, \Kapp}  -\metric{\nabla_s\Kapp,\nabla_s\phi} \metric{\psi, \Kapp}\\
&\qquad  -\frac{3}{2}|\Kapp|^2\metric{\nabla_s \psi, \nabla_s \phi}+ \metric{\phi,\Kapp}\metric{\Kapp, \nabla_s^2 \psi}+ \metric{\nabla_s \psi, \Kapp}\metric{\Kapp,\nabla_s \phi}\\
&\qquad+\metric{ \nabla_s^2 \phi,\Kapp}\metric{\Kapp,\psi} +|\Kapp|^2 \metric{\Kapp, \psi} \metric{\Kapp,\phi} \diff s
         \end{split}
\end{align}
which yields the claim after rearranging the terms. 
\end{proof}


 \subsection{Calculation of the second variation of the length functional }

\begin{proof}[Proof of Corollary \ref{cor:2varl}]\label{proof:2ndvarlambda}
For simplicity of notation in this proof we denote the Euclidean scalar-product simply by $\metric{\cdot,\cdot}$ and we write $s=s_{\bar f}$ and $\Kapp=\Kapp_{\bar{f}}$. 
 The length of $\bar{f}+\eps \phi$ is given by 
$$ \mathcal{L}(\bar{f}+\eps \phi) = \int_I |\partial_x \bar{f} + \eps \partial_x \phi| \diff x $$
and hence by \eqref{eq:varnorm} $\mathcal{L}'(\bar{f}) \phi = -\int_I \metric{\Kapp,\phi} \diff s $. Thus the $L^2$-gradient of the length functional is given by $\nabla_{L^2(ds)}\mathcal{L}(\bar{f}) = - \Kapp$.
Furthermore, for the second variation we find by \eqref{eq:varform},\eqref{eq:varkapp} and orthogonality that
\begin{align*}
\mathcal{L}''(\bar{f})[\psi, \phi] &= - 
\left. \frac{d}{d\varepsilon} \right|_{\varepsilon=0} 
\int_I  \metric{\Kapp_{\bar{f}+\eps \phi},\psi} \diff s_{\bar{f}+\eps \phi}  
\\
& = -\int_I (\Kapp \cdot \psi) (-\phi \cdot \Kapp) + \metric{ \nabla_s^2 \phi - (\nabla_s \phi\cdot\Kapp) \partial_s \bar{f} + (\phi \cdot \Kapp)\Kapp, \psi} \diff s\\
&= -\int\metric{\nabla_s^2 \phi, \psi} 
\diff s = \int \metric{\nabla_s \phi, \nabla_s \psi}  \diff s,
\end{align*}
applying \eqref{eq:partintnabla} with $\eta = \nabla_s \phi$ and $\xi = \psi$ and the facts that $\psi$ is in $H^1_0$ and $\nabla_s \phi$, $\psi$ are  orthogonal to $\partial_s \bar{f}$.
Together with \eqref{eq:2ndvar} this shows the claim.
\end{proof}

\subsection{Technical lemmas}
\begin{lemma}
Let $f:I \to \R^d$ be a smooth regular curve, $\phi:I \to \R^d$ be a smooth vector field, and $s=s_{f}$. Then
\begin{align}\label{eqa1}
\partial_s \phi & = \nabla_s \phi + \langle \partial_s \phi, \partial_s f \rangle_\text{euc} \partial_s f  \, ,\\ \label{eqa2}
\partial_s^2 \phi & = \nabla_s^2 \phi +  \langle \partial_s^2  \phi, \partial_s f \rangle_\text{euc} \partial_s f + \langle \partial_s \phi, \partial_s f \rangle_\text{euc} \Kapp \, .
\end{align}
If $\phi$ is normal to $f$, i.e. $\langle  \phi, \partial_s f \rangle_\text{euc}=0$  on $I$, then
\begin{align}\label{eqa3}
\partial_s \phi & = \nabla_s \phi - \langle  \phi, \Kapp \rangle_\text{euc} \partial_s f  \, ,\\ \label{eqa4}
\partial_s^2 \phi & = \nabla_s^2 \phi -2 \langle \nabla_s \phi, \Kapp \rangle_\text{euc} \partial_s f 
 - \langle \phi, \nabla_s \Kapp \rangle_\text{euc} \partial_s f   - \langle \phi, \Kapp \rangle_\text{euc} \Kapp \, ,
\end{align}
and 
\begin{align}\label{eqa5}
\partial_s^3 \phi & =    \nabla_s^3 \phi -3\langle \nabla_s^2 \phi, \Kapp \rangle_\text{euc} \partial_s f -3\langle \nabla_s \phi, \nabla_s \Kapp \rangle_\text{euc} \partial_s f -\langle \phi, \nabla_s^2 \Kapp \rangle_\text{euc} \partial_s f \\ \nonumber
& \quad -3 \langle \nabla_s \phi, \Kapp \rangle_\text{euc} \Kapp  -2 \langle \phi, \nabla_s \Kapp \rangle_\text{euc} \Kapp - \langle \phi, \Kapp \rangle_\text{euc} \partial_s \Kapp \, ,\\ \label{eqa6}
\partial_s^4 \phi & =  \nabla_s^4 \phi  - 4 \langle  \nabla_s^3 \phi, \Kapp \rangle_\text{euc} \partial_s f  - 6 \langle  \nabla_s^2 \phi, \nabla_s \Kapp \rangle_\text{euc} \partial_s f - 6 \langle  \nabla_s^2 \phi, \Kapp \rangle_\text{euc} \Kapp \\
\nonumber
& \quad   - 4 \langle  \nabla_s \phi, \nabla_s^2 \Kapp \rangle_\text{euc} \partial_s f   - 8 \langle  \nabla_s \phi, \nabla_s \Kapp \rangle_\text{euc} \Kapp    - 4 \langle  \nabla_s \phi, \Kapp \rangle_\text{euc} \partial_s \Kapp\\
\nonumber
& \quad    - \langle  \phi, \nabla_s^3 \Kapp \rangle_\text{euc} \partial_s f  - 3 \langle  \phi, \nabla_s^2 \Kapp \rangle_\text{euc} \Kapp   -3 \langle  \phi, \nabla_s \Kapp \rangle_\text{euc} \partial_s \Kapp     - \langle  \phi, \Kapp \rangle_\text{euc}\partial_s^2 \Kapp\, .
\end{align}
\end{lemma}
\begin{proof}
The first formula is exactly the definition of the operator $\nabla_s$. Differentiating once more we obtain
\begin{align*}
\partial_s^2 \phi & = \nabla_s^2 \phi + \langle \partial_s \nabla_s \phi, \partial_s f \rangle_\text{euc} \partial_s f + \langle \partial_s^2 \phi, \partial_s f \rangle_\text{euc} \partial_s f \\
& \quad + \langle \partial_s \phi, \Kapp \rangle_\text{euc} \partial_s f + \langle \partial_s \phi, \partial_s f \rangle_\text{euc} \Kapp \\
& = \nabla_s^2 \phi +  \langle \partial_s^2  \phi, \partial_s f \rangle_\text{euc} \partial_s f + \langle \partial_s \phi, \partial_s f \rangle_\text{euc} \Kapp, 
\end{align*} 
since $\partial_s^2 f =\Kapp$.

If $\phi$ is normal to $f$, then $ \langle \partial_s \phi, \partial_s f \rangle_\text{euc} = -  \langle \phi, \Kapp \rangle_\text{euc}$ and hence \eqref{eqa3} follows directly from \eqref{eqa1}. Differentiating \eqref{eqa3} once we get
\begin{align*}
\partial_s^2 \phi & = \nabla_s^2 \phi + \langle  \partial_s \nabla_s \phi, \partial_s f \rangle_\text{euc} \partial_s f - \langle  \nabla_s \phi, \Kapp \rangle_\text{euc} \partial_s f \\
& \quad  - \langle  \phi, \nabla_s \Kapp \rangle_\text{euc} \partial_s f - \langle  \phi, \Kapp \rangle_\text{euc} \Kapp \\
 & = \nabla_s^2 \phi - 2 \langle  \nabla_s \phi, \Kapp \rangle_\text{euc} \partial_s f   - \langle  \phi, \nabla_s \Kapp \rangle_\text{euc} \partial_s f - \langle  \phi, \Kapp \rangle_\text{euc} \Kapp \, ,
\end{align*}
that is \eqref{eqa4}. Proceeding similarly,
\begin{align*}
\partial_s^3 \phi & = \partial_s \nabla_s^2 \phi 
- 2 \langle  \nabla_s^2 \phi, \Kapp \rangle_\text{euc} \partial_s f  - 2 \langle  \nabla_s \phi, \nabla_s \Kapp \rangle_\text{euc} \partial_s f  - 2 \langle  \nabla_s \phi, \Kapp \rangle_\text{euc} \Kapp \\
& \quad   - \langle \nabla_s \phi, \nabla_s \Kapp \rangle_\text{euc} \partial_s f  - \langle  \phi, \nabla_s^2 \Kapp \rangle_\text{euc} \partial_s f 
 - \langle  \phi, \nabla_s \Kapp \rangle_\text{euc} \Kapp \\
& \quad  - \langle \nabla_s \phi, \Kapp \rangle_\text{euc} \Kapp - \langle  \phi,\nabla_s \Kapp \rangle_\text{euc} \Kapp- \langle  \phi, \Kapp \rangle_\text{euc}\partial_s \Kapp \\
&= \nabla_s^3 \phi  - 3 \langle  \nabla_s^2 \phi, \Kapp \rangle_\text{euc} \partial_s f  - 3 \langle  \nabla_s \phi, \nabla_s \Kapp \rangle_\text{euc} \partial_s f  \\
& \quad - 3 \langle  \nabla_s \phi, \Kapp \rangle_\text{euc} \Kapp  - \langle  \phi, \nabla_s^2 \Kapp \rangle_\text{euc} \partial_s f 
 - 2 \langle  \phi, \nabla_s \Kapp \rangle_\text{euc} \Kapp   - \langle  \phi, \Kapp \rangle_\text{euc}\partial_s \Kapp \, ,
\end{align*}
and
\begin{align*}
\partial_s^4 \phi & =\partial_s \nabla_s^3 \phi  - 3 \langle  \nabla_s^3 \phi, \Kapp \rangle_\text{euc} \partial_s f  - 3 \langle  \nabla_s^2 \phi, \nabla_s \Kapp \rangle_\text{euc} \partial_s f - 3 \langle  \nabla_s^2 \phi, \Kapp \rangle_\text{euc} \Kapp \\
& \quad  - 3 \langle  \nabla_s^2 \phi, \nabla_s \Kapp \rangle_\text{euc} \partial_s f   - 3 \langle  \nabla_s \phi, \nabla_s^2 \Kapp \rangle_\text{euc} \partial_s f   - 3 \langle  \nabla_s \phi, \nabla_s \Kapp \rangle_\text{euc} \Kapp  \\
& \quad   - 3 \langle  \nabla_s^2 \phi, \Kapp \rangle_\text{euc} \Kapp - 3 \langle  \nabla_s \phi, \nabla_s \Kapp \rangle_\text{euc} \Kapp - 3 \langle  \nabla_s \phi, \Kapp \rangle_\text{euc} \partial_s \Kapp\\
& \quad   - \langle \nabla_s \phi, \nabla_s^2 \Kapp \rangle_\text{euc} \partial_s f  - \langle  \phi, \nabla_s^3 \Kapp \rangle_\text{euc} \partial_s f  - \langle  \phi, \nabla_s^2 \Kapp \rangle_\text{euc} \Kapp \\
& \quad  - 2 \langle \nabla_s \phi, \nabla_s \Kapp \rangle_\text{euc} \Kapp - 2 \langle  \phi, \nabla_s^2 \Kapp \rangle_\text{euc} \Kapp - 2 \langle  \phi, \nabla_s \Kapp \rangle_\text{euc} \partial_s \Kapp\\
& \quad - \langle \nabla_s \phi, \Kapp \rangle_\text{euc}\partial_s \Kapp - \langle  \phi, \nabla_s \Kapp \rangle_\text{euc}\partial_s \Kapp
   - \langle  \phi, \Kapp \rangle_\text{euc}\partial_s^2 \Kapp\\
& =\nabla_s^4 \phi  - 4 \langle  \nabla_s^3 \phi, \Kapp \rangle_\text{euc} \partial_s f  - 6 \langle  \nabla_s^2 \phi, \nabla_s \Kapp \rangle_\text{euc} \partial_s f - 6 \langle  \nabla_s^2 \phi, \Kapp \rangle_\text{euc} \Kapp \\
& \quad   - 4 \langle  \nabla_s \phi, \nabla_s^2 \Kapp \rangle_\text{euc} \partial_s f   - 8 \langle  \nabla_s \phi, \nabla_s \Kapp \rangle_\text{euc} \Kapp    - 4 \langle  \nabla_s \phi, \Kapp \rangle_\text{euc} \partial_s \Kapp\\
& \quad    - \langle  \phi, \nabla_s^3 \Kapp \rangle_\text{euc} \partial_s f  - 3 \langle  \phi, \nabla_s^2 \Kapp \rangle_\text{euc} \Kapp   -3 \langle  \phi, \nabla_s \Kapp \rangle_\text{euc} \partial_s \Kapp     - \langle  \phi, \Kapp \rangle_\text{euc}\partial_s^2 \Kapp \, .
\end{align*}
\end{proof}

\subsection{Completion of the proof of Lemma~\ref{lemma:reparametrization}}
\label{app:proofs}
\begin{proof}[Second part of the claim of Lemma \ref{lemma:reparametrization}]
For the second part of the claim we need to prove some estimates. We may rewrite \eqref{auf1} as
\begin{align*}
\tilde{\varphi}'  &  =\frac{ |\partial_x \bar{f}|^2 - \metric[\text{euc}]{(\partial_x\bar{f})\circ (\id_I + \tilde\varphi) , \partial_x \bar{f}} -  \metric[\text{euc}]{\psi\circ (\id_I + \tilde\varphi) , \partial_x^2 \bar{f}}}{\metric[\text{euc}]{\partial_x(\bar{f}+\psi)\circ (\id_I + \tilde\varphi) , \partial_x \bar{f}}} \\
& \quad + \frac{  -  \metric[\text{euc}]{\bar{f}\circ (\id_I + \tilde\varphi) - \bar{f}, \partial_x^2 \bar{f}}- \metric[\text{euc}]{(\partial_x\psi)\circ (\id_I + \tilde\varphi) , \partial_x \bar{f}}}{\metric[\text{euc}]{\partial_x(\bar{f}+\psi)\circ (\id_I + \tilde\varphi) , \partial_x \bar{f}}}\, .
\end{align*}
In order to proceed we recall that if $g_1,g_2 \in H^3$ then $\| g_1 g_2\|_{H^3} \leq c \| g_1\|_{H^3}\|g_2\|_{H^3}$. If furthermore $g_2 \geq \frac 12 |\partial_x \bar{f}|^2$, then
$$ \Big \|\frac{1}{g_2} \Big \|_{H^3} \leq C(\bar{f})  \left(1+ \|g_2\|_{H^3} + \|g_2\|^2_{H^3} +\|g_2\|^3_{H^3} \right)  \leq C(\bar{f})  \left(1+ \|g_2\|^3_{H^3} \right) \, .$$
While considering the composition if $\frac12 \leq |g_2'|\leq \frac32$, then
$$ \|g_1 \circ g_2\|_{H^3} \leq C  \|g_1\|_{H^3} (1+ \|g_2\|_{H^3})\, .$$ 
By using several time these estimates we find writing 
\begin{align*}
\partial_x \bar{f}\circ (\id_I + \tilde\varphi)  -\partial_x \bar{f} =  
\int_{\id_I}^{\id_I+\tilde\varphi}  \partial_x^2 \bar{f}(y) \; \diff y
\end{align*}
that
\begin{align*}
\| \partial_x \bar{f}\circ (\id_I + \tilde\varphi)  -\partial_x \bar{f}  \|_{H^3} & \leq  C(\bar{f}) (\| \tilde\varphi\|_{H^3} +\|\tilde \varphi\|_{L^2}^\frac{1}{2} ),\\
\| |\partial_x \bar{f}|^2 - \metric[\text{euc}]{\partial_x(\bar{f})\circ (\id_I + \tilde\varphi) , \partial_x \bar{f}} \|_{H^3} & \leq  C(\bar{f})  (\| \tilde\varphi\|_{H^3} +\|\tilde \varphi\|_{L^2}^\frac{1}{2} ),\\
\| \bar{f}\circ (\id_I + \tilde\varphi)  -\bar{f}  \|_{H^3} & \leq  C(\bar{f}) \| \tilde\varphi\|_{H^3},\\
\| \metric[\text{euc}]{\bar{f}\circ (\id_I + \tilde\varphi) - \bar{f}, \partial_x^2 \bar{f}}\|_{H^3} & \leq  C(\bar{f}) \| \tilde\varphi\|_{H^3}, \\
\|\metric[\text{euc}]{\psi\circ (\id_I + \tilde\varphi) , \partial_x^2 \bar{f}} \|_{H^3}  &  \leq C(\bar{f}) \| \psi\|_{H^3} (1 + \| \tilde\varphi\|_{H^3}),\\
\|\metric[\text{euc}]{\partial_x \psi\circ (\id_I + \tilde\varphi) , \partial_x \bar{f}} \|_{H^3}  &  \leq C(\bar{f}) \| \psi\|_{H^4} (1+  \|\tilde\varphi\|_{H^3})  ,
\end{align*}
and
\begin{align*}
& \qquad \Big \| \frac{ 1}{\metric[\text{euc}]{\partial_x(\bar{f}+\psi)\circ (\id_I + \tilde\varphi) , \partial_x \bar{f}}}  \Big \|_{H^3} \\
& \leq C(\bar{f}) (1+ \| \metric[\text{euc}]{\partial_x \bar{f} \circ (\id_I + \tilde\varphi) , \partial_x \bar{f}} +\metric[\text{euc}]{ (\partial_x\psi)\circ (\id_I + \tilde\varphi) , \partial_x \bar{f}} \|_{H^3}^3)\\
& \leq C(\bar{f}) (1+  \| \tilde\varphi\|_{H^3}^3 + \| \psi\|_{H^4}^3 (1+  \| \tilde\varphi\|_{H^3})^3) \, .
\end{align*}
From these inequalities and the fact that $\| \tilde{\varphi}\|_{H^3}$ depends continuously on $\| \psi\|_{H^4}$ we get
$$\| \tilde{\varphi}\|_{H^4} \leq C(\bar{f}, \| \psi\|_{H^4}) \, ,$$
with a constant depending continuously on $\| \psi\|_{H^4}$ and such that when $\| \psi\|_{H^4} \searrow 0$ then also 
$\| \tilde{\varphi}\|_{H^4} \searrow 0$. Since $\Phi = \id+\tilde\varphi$ and $\phi = (\bar{f}+\psi)\circ \Phi - \bar{f}$ with estimates similar to the ones above, we find
\[ \| \phi \|_{H^4} \leq C(\bar{f}) \|\tilde\varphi\|_{H^4} + C(\bar{f})\| \psi\|_{H^4}(1+  \| \tilde\varphi\|_{H^4}^{2}) \, .\]
This gives the claim. 
\end{proof}

\bibliography{Lib}
\bibliographystyle{alpha}


\begin{small}

\noindent \textit{Anna Dall'Acqua, Adrian Spener},

 \noindent Universit\"at Ulm, Helmholtzstra\ss e 18, 89081 Ulm, Germany
 
 \noindent
\texttt{anna.dallacqua@uni-ulm.de, adrian.spener@uni-ulm.de}
 
\smallskip

\noindent \textit{Paola Pozzi}, 

\noindent Universit\"at Duisburg-Essen, Mathematikcarr\'ee, 
Thea-Leymann-Stra\ss e 9, 45127 Essen, Germany

\noindent \texttt{paola.pozzi@uni-due.de}
\end{small}

\end{document}